\documentclass[11pt,twoside]{amsart}
\usepackage{latexsym,amssymb,amsmath}

\numberwithin{equation}{section}
\newtheorem{theorem}{Theorem}[section]
\newtheorem{lemma}[theorem]{Lemma}
\newtheorem{proposition}[theorem]{Proposition}
\newtheorem{corollary}[theorem]{Corollary}
\newtheorem{conjecture}[theorem]{Conjecture}

\theoremstyle{definition}
\newtheorem{definition}[theorem]{Definition} 
\newtheorem{remark}[theorem]{Remark}
\newtheorem{example}[theorem]{Example}

\begin{document}

\title[Cohen-Macaulay, Shellable and unmixed clutters]
{Cohen-Macaulay, Shellable and unmixed clutters 
with a perfect matching of K\"onig type}
 
\author{Susan Morey}
\address{Department of Mathematics \\
Texas State University\\
601 University Drive\\ 
San Marcos, TX 78666.}
\email{morey@txstate.edu}

\author{Enrique Reyes}
\address{Departamento de Ciencias B\'asicas, Unidad Profesional
Interdisciplinaria en Ingenier\'\i a y Tecnologias Avanzadas del IPN, 
UPIITA, Av. IPN 2580, Col. Barrio la Laguna Ticom\'an, 07340 M\'exico
City, D.F.} 
\email{ereyes@math.cinvestav.mx}
\author{Rafael H. Villarreal.}
\address{
Departamento de
Matem\'aticas\\
Centro de Investigaci\'on y de Estudios
Avanzados del
IPN\\
Apartado Postal
14--740 \\
07000 Mexico City, D.F.
}
\email{vila@math.cinvestav.mx}

  
\begin{abstract} Let $\mathcal{C}$ be a clutter with a perfect
matching $e_1,\ldots,e_g$ of K\"onig type and let 
$\Delta_\mathcal{C}$ be the Stanley-Reisner complex of the edge 
ideal of $\mathcal{C}$. If all
c-minors of $\mathcal{C}$ have a free vertex and $\mathcal{C}$ is
unmixed, we show that $\Delta_\mathcal{C}$ is pure shellable. We are
able 
to describe, in combinatorial and algebraic terms, when
$\Delta_\mathcal{C}$ is pure. 
If $\mathcal{C}$ has no cycles of length $3$ or $4$, then it is shown 
that $\Delta_\mathcal{C}$ is pure if and only if $\Delta_\mathcal{C}$
is pure shellable (in this case $e_i$ has a free vertex for all
$i$), and that $\Delta_\mathcal{C}$ is pure if and only 
if for any two edges   
$f_1,f_2$ of $\mathcal{C}$ and 
for any $e_i$, one has that $f_1\cap e_i\subset f_2\cap e_i$ or
$f_2\cap e_i\subset f_1\cap e_i$. It is also shown that this ordering
condition implies that $\Delta_\mathcal{C}$ is pure shellable, 
without any assumption on the cycles of $\mathcal{C}$. 
Then we prove that complete 
admissible uniform clutters and their Alexander duals are unmixed. 
In addition, the edge ideals of complete admissible uniform clutters 
are facet ideals of shellable simplicial
complexes, they are Cohen-Macaulay, and they have linear resolutions. 
Furthermore if $ \mathcal{C}$ is admissible and
complete, then $\mathcal{C}$ is unmixed. 
We characterize certain conditions that occur in a Cohen-Macaulay
criterion for 
bipartite graphs of Herzog and Hibi, and extend some results of
Faridi---on the structure of unmixed simplicial trees---to clutters
with the K\"onig property without $3$-cycles or $4$-cycles. 
\end{abstract}
 
\maketitle

\section{Introduction}

A {\it clutter\/} $\mathcal{C}$ with finite vertex set $X$ is a 
family of subsets of $X$, 
called edges, none of which is included in another. The set of
vertices and edges of $\mathcal{C}$ are denoted by $V(\mathcal{C})$
and $E(\mathcal{C})$ respectively. Clutters are special types of
hypergraphs. The set 
of edges of a clutter can be viewed as the set of facets of a
simplicial complex. A basic example  
of a clutter is a graph. For a thorough study of clutters and
hypergraphs from the point of view of combinatorial optimization see
\cite{cornu-book,Schr2}.  

Let $\mathcal{C}$ be a clutter with finite vertex set 
$X=\{x_1,\ldots,x_n\}$. We shall always assume that 
$\mathcal{C}$ has no isolated vertices, i.e., each vertex occurs in
at least one edge. Let $R=K[x_1,\ldots,x_n]$ be a polynomial ring 
over a field $K$. 
The {\it edge ideal\/} of $\mathcal{C}$, 
denoted by $I(\mathcal{C})$, is the ideal of $R$
generated by all monomials $\prod_{x_i\in e}x_i=x_e$ such 
that $e\in E(\mathcal{C})$. The assignment $\mathcal{C}\mapsto
I(\mathcal{C})$ 
establishes a natural one to one
correspondence between the family of clutters and the family of 
square-free monomial ideals. Edge ideals of clutters
are also called {\it facet ideals} \cite{Faridi}. A subset $F$ of $X$
is called {\it independent\/} or {\it stable\/} if $e\not\subset F$
for any  
$e\in E(\mathcal{C})$. The dual concept of an independent vertex set
is a 
{\it vertex cover\/}, i.e., a subset $C$ of $X$ is a vertex cover of
$\mathcal{C}$  
if and only if $X\setminus C$ is an independent vertex set. The
number of vertices in a minimum vertex cover of $\mathcal{C}$ 
is called the {\it covering
number\/} of $\mathcal{C}$, and this number coincides with ${\rm
ht}\, I(\mathcal{C})$, the {\it height\/} of $I(\mathcal{C})$.
The Stanley-Reisner complex of
$I(\mathcal{C})$, denoted by $\Delta_\mathcal{C}$, is the simplicial
complex whose faces are the independent vertex sets of $\mathcal{C}$.
Recall 
that $\Delta_\mathcal{C}$ is called {\it pure\/} 
if all maximal independent vertex sets of $\mathcal{C}$, with respect
to inclusion, have the same number of
elements. If $\Delta_\mathcal{C}$ is pure (resp. Cohen-Macaulay,
Shellable), we say that $\mathcal{C}$ is
{\it unmixed\/} (resp. Cohen-Macaulay, Shellable). A clutter has the
{\it K\"onig property\/} if the maximum number of pairwise disjoint
edges equals the covering number. 
A {\it perfect matching of $\mathcal{C}$ of K\"onig type\/} is a
collection $e_1,\ldots,e_g$ of pairwise disjoint 
edges whose union is $X$ and
such that $g$ is the height of $I(\mathcal{C})$. 
Any unmixed clutter
with the K\"onig property and without isolated vertices has a perfect
matching of K\"onig type (Lemma~\ref{um+konig-kt}). 

We are 
interested in determining what families of clutters have the property
that $\Delta_\mathcal{C}$ is pure, Cohen-Macaulay, 
or Shellable in the non-pure sense of Bj\"orner-Wachs \cite{BW}. 
The last two properties have been extensively studied, see
\cite{BHer,ITG,Stanley,monalg} and the references there, but to the best
of our knowledge the first property
has not been studied much except for the case of graphs
\cite{plummer-unmixed,plummer-survey,ravindra,unmixed}. 
The aim of
this paper is to examine these three properties when $\mathcal{C}$ 
has a perfect matching of K\"onig type or when $\mathcal{C}$ has the 
K\"onig property.

The contents of this paper are as follows. Let $\mathcal{C}$ be a
clutter with a 
perfect matching $e_1,\ldots,e_g$ of K\"onig type and let
$I(\mathcal{C})$ be its edge ideal.  The main theorem in Section 2 is
a combinatorial description of the unmixed property of $\mathcal{C}$, 
along with an equivalent algebraic formulation. Before stating the 
theorem, recall that the {\it support\/} of
$x^a=x_1^{a_1}\cdots x_n^{a_n}$, denoted by 
${\rm supp}(x^a)$, is the set of $x_i$ such that $a_i>0$. The
{\it colon ideal\/} $(x^a\colon x^b)$ is the set of  $f$ in $R$ such 
that $fx^b$ is in $(x^a)$. The {\it colon ideal\/} 
$(I(\mathcal{C})^2\colon x_{e_i})$ is defined similarly.  

\noindent {\bf Theorem~\ref{unmixed-clutter-pm}.} The following 
conditions are equivalent:
\begin{itemize}
\item[(a)] $\mathcal{C}$ is unmixed.
\item[(b)] For any two edges $e \neq e'$ and for any two distinct
vertices  
$x\in e$, $y\in e'$ contained in some $e_i$, one has that 
$(e\setminus\{x\})\cup (e'\setminus\{y\})$ contains an edge. 
\item[(c)] For any two edges $e\neq e'$ and for any $T\subset e_i$
such that $x_T$ divides $x_ex_{e'}$, one has that 
${\rm supp}(x_ex_{e'}/x_T)$ contains an edge. 
\item[(d)] For any two edges $e\neq e'$ and for any $e_i$,
$(x_ex_{e'}\colon x_{e_i}) 
\subset I(\mathcal{C})$.
\item[(e)] $I(\mathcal{C})=(I(\mathcal{C})^2\colon
x_{e_1})+\cdots+(I(\mathcal{C})^2\colon x_{e_g})$.
\end{itemize}

\noindent This generalizes to balanced
clutters  (see Definition~\ref{bala-clu}) and beyond an unmixedness 
criterion of \cite{unmixed} valid only for bipartite graphs 
(Corollary~\ref{char-bal-um}). 

The notions of minor and c-minor play a prominent role in 
combinatorial optimization \cite{cornu-book}. The precise definitions
of these notions can be found in Section~\ref{scpm}. Roughly speaking
a minor (c-minor) is obtained  
from $I(\mathcal{C})$ by making any sequence of variables equal to
$1$ or $0$ (resp. equal to $1$ only). From the algebraic point of
view, a c-minor corresponds to a colon operation or 
localization of $I(\mathcal{C})$. In Theorem~\ref{c-minors} we show
that for a clutter with a perfect matching of K\"{o}nig type, if all
c-minors of $\mathcal{C}$ have a free  
vertex, i.e., a vertex that occurs in one edge only, 
and $\mathcal{C}$ is unmixed, then $\Delta_\mathcal{C}$
is pure shellable. This complements a result of 
\cite{bipartite-scm} showing that if all minors of an arbitrary
clutter $\mathcal{C}$ have a free vertex,  then $\Delta_\mathcal{C}$
is shellable. Using this free vertex property, we show in
Theorem~\ref{ordering+matching-shellable} that if for any two edges   
$f_1,f_2$ of  $\mathcal{C}$ and  for any $e_i$, one has that 
$f_1\cap e_i\subset f_2\cap e_i$ or $f_2\cap e_i\subset f_1\cap e_i$,
then  $\Delta_{\mathcal C}$ is pure shellable. Note that this
ordering property on the edges implies $\mathcal{C}$ is unmixed, as
is seen in Theorem~\ref{jun1-07-1}.    

An additional property is needed to guarantee that an unmixed clutter
will have the above ordering property. Let $A=(a_{ij})$ be the 
{\it incidence\/} matrix
of $\mathcal{C}$. Recall that $a_{ij}=1$ if $x_i\in g_j$ and
$a_{ij}=0$ otherwise, where $g_1,\ldots,g_q$ are the edges of
$\mathcal{C}$. In Theorem~\ref{jun1-07} we assume that $\mathcal{C}$
has no cycles of length $3$ or $4$, i.e., $A$ has no square  
submatrix of order $3$ or $4$ with exactly two $1$'s in 
each row and column, and then show that if $\mathcal{C}$ is unmixed,
then  
for any two edges $f_1,f_2$ of $\mathcal{C}$ and for any $e_i$, one
has that  
$f_1\cap e_i\subset f_2\cap e_i$ or $f_2\cap e_i\subset f_1\cap e_i$.  
This ordering property was shown to hold for the clutter of 
facets of any unmixed simplicial tree \cite[Remark~7.2,
Corollary~7.8]{faridijct}. Thus our result  
is a wide generalization of this fact because simplicial trees 
are acyclic clutters \cite{hhtz}. In addition, when $\mathcal{C}$ is
unmixed and has no cycles of length three or four, we show that
$\Delta_{\mathcal C}$ is pure shellable (Theorem~\ref{jun2-07}) and that  
$e_i$ has a free vertex for all $i$
(Proposition~\ref{structure-theorem}). Then we give a far reaching
generalization of Faridi's characterization of unmixed 
simplicial trees \cite{faridijct} (see Corollary~\ref{char-tbc}) and show
some applications of these results to totally balanced clutters 
(Corollary~\ref{structureoftrees}).

In Section~\ref{clutters-pm} we introduce the notion of an admissible
clutter. The notion of an admissible clutter was inspired by a certain 
ordering condition that occurs in a Cohen-Macaulay criterion for
bipartite graphs of Herzog and Hibi \cite{herzog-hibi-crit} (see
condition {\rm 
($\mathrm{h}_1$)} below). We show that any complete admissible 
clutter is unmixed (Proposition~\ref{complete-case-enrique}) and that 
the edge ideal of any complete admissible uniform clutter is the facet
ideal of a shellable complex (Theorem~\ref{susan-shellable}). 
A clutter is called {\it uniform\/} if all its edges have the same
size. It is 
shown in Lemma~\ref{susan-duality} that complete admissible uniform
clutters are closed under 
taking Alexander duals. This allows us to prove
Theorem~\ref{complete-CM}: If $\mathcal{C}$ is a complete admissible
uniform clutter, then $R/I({\mathcal C})$ is Cohen-Macaulay and has
a linear resolution.    
An interesting problem that remains
unsolved is whether an unmixed admissible clutter is Cohen-Macaulay 
(Conjecture~\ref{problem}). For bipartite graphs this problem has a
positive answer (Theorem~\ref{herzog-hibi}, \cite{herzog-hibi-crit}).

Section~\ref{pm-bipartite-konig} is devoted to bipartite graphs with
a perfect matching of K\"onig type. 
An unmixed bipartite graph without isolated vertices will always have
this type of matching by K\"onig's theorem \cite{Schr2}. 
Bipartite Cohen-Macaulay graphs have been studied in
\cite{carra-ferrarello,EV,herzog-hibi-crit,monalg}. In \cite{EV} it is
shown that $G$ is a 
Cohen-Macaulay graph if and only if $\Delta_G$ is pure 
shellable. In \cite{bipartite-scm} a classification of all 
sequentially Cohen-Macaulay bipartite graphs is given. In particular,
it is 
shown that $\Delta_G$ is shellable if and only if $R/I(G)$ is
sequentially Cohen-Macaulay.   

Let $G$ be a bipartite graph and let $V_1=\{x_1,\ldots,x_g\}$ and 
$V_2=\{y_1,\ldots,y_g\}$ be a bipartition of $G$ such 
that $\{x_i,y_i\}\in E(G)$ for all $i$. We examine the conditions
($\mathrm{h}_1$): ``if $\{x_i, y_j\} \in E(G)$, then $i \leq j$", and 
($\mathrm{h}_2$): ``if $\{x_i, y_j\}$ and $\{x_j, y_k\}$ are in
$E(G)$ and $i<j<k$, then $\{x_i, y_k\}\in E(G)$" that occur in the 
Herzog and Hibi criterion for Cohen-Macaulay bipartite graphs
\cite{herzog-hibi-crit}. See Theorem~\ref{herzog-hibi} for a precise 
statement of this criterion. Some characterizations of these  
conditions have been shown by Yassemi (personal communication), and
by  Carr\`a Ferro and Ferrarello \cite{carra-ferrarello}. These 
conditions have also been examined in \cite{bipartite-scm} from the
point of view of digraphs following ideas introduced in
\cite{carra-ferrarello}. Our main result of
Section~\ref{pm-bipartite-konig} shows that 
condition ($\mathrm{h}_1$) holds if and only if the subcomplex
generated  
by the facets of maximum dimension of $\Delta_G$ is shellable
(Theorem~\ref{ordering-char}). We recover a result of \cite{unmixed} 
describing all unmixed bipartite graphs in combinatorial terms
(Corollary~\ref{unmixed-bip}). In particular it follows that in the
Herzog 
and Hibi criterion (Theorem~\ref{herzog-hibi}) 
we can replace condition ($\mathrm{h}_2$) by
condition ($\mathrm{h}_2'$): ``$G$ is unmixed''. In
Corollary~\ref{herzog-hibi-ref} we give a variation of this 
criterion. 

The natural 
generalization of a bipartite graph is a balanced clutter, i.e., a
clutter without odd cycles. 
It turns out that the ordering criterion that Herzog and Hibi used to
classify Cohen-Macaulay bipartite graphs does not extend to 
Cohen-Macaulay balanced 
clutters (Example~\ref{counterexample}). 

\section{Shellable clutters with a perfect matching}\label{scpm} 
Let $\mathcal{C}$ be a clutter on the vertex set
$X=\{x_1,\ldots,x_n\}$ and let $I=I(\mathcal{C})$ be its edge ideal.
A {\it contraction\/} (resp. {\it deletion}) of $I$ is an ideal of
the form  
$(I\colon x_i)$ (resp. $J=I\cap
K[x_1,\ldots,\widehat{x}_i,\ldots,x_n]$) for some $x_i$, where 
$(I\colon x_i):=\{f\in R\vert\, f x_i\in I\}$ is the standard colon
operation in ideal theory. The ideal 
$I$ is regarded as a contraction. The clutter associated to 
the square-free monomial ideal $(I\colon x_i)$ (resp. $J$) is denoted
by $\mathcal{C}/x_i$  
(resp. $\mathcal{C}\setminus x_i$). A {\it c-minor\/} 
(resp. {\it d-minor\/}) of $I$ is an ideal
obtained from $I$ by a sequence of contractions (resp. deletions). If
a c-minor $I'$ contains a variable $x_i$ and we remove this variable
from $I'$, we still consider the new ideal a c-minor of $I$. A {\it
minor\/} of $I$
is an ideal obtained from $I$ by a sequence of deletions 
and contractions in any order. A {\it minor\/} (resp. c-minor) of
$\mathcal{C}$ is any  
clutter that correspond to a minor (resp. c-minor) of $I$. This
terminology is 
consistent with that of \cite[p.~23]{cornu-book}. A vertex $x$ of
$\mathcal{C}$ is called {\it isolated\/} if $x$ does not occur in any
edge of $\mathcal{C}$. A subset $C\subset X$ is a 
{\it minimal vertex cover\/} of the clutter $\mathcal{C}$ if: 
($\mathrm{c}_1$) every edge of $\mathcal{C}$ contains at least one
vertex of $C$,  
and ($\mathrm{c}_2$) there is no proper subset of $C$ with the first 
property. If $C$ only satisfies condition ($\mathrm{c}_1$), then $C$ is 
called a {\it vertex cover\/} of $\mathcal{C}$. Recall that 
$\mathfrak{p}$ is a minimal prime of $I =I(\mathcal{C})$ if and only if 
$\mathfrak{p}=(C)$ for some minimal vertex cover $C$ of $\mathcal{C}$ 
\cite[Proposition~6.1.16]{monalg}. Thus the primary decomposition of
the edge ideal of $\mathcal{C}$ is given by
$$
I(\mathcal{C})=(C_1)\cap (C_2)\cap\cdots\cap (C_p),
$$
where $C_1,\ldots,C_p$ are the minimal vertex 
covers of $\mathcal{C}$. In particular observe that the height of 
$I(\mathcal{C})$ equals the number of vertices in a minimum vertex
cover of $\mathcal{C}$. Note that the facets of $\Delta_\mathcal{C}$ 
are $X\setminus C_1,\ldots,X\setminus C_p$. Thus $\mathcal{C}$ is
{\it unmixed}, equivalently $\Delta_\mathcal{C}$ is pure, if
and only if all minimal vertex covers of $\mathcal{C}$ have the same
size.  

\begin{definition} A {\it perfect matching of K\"onig type\/} of
$\mathcal{C}$ is a collection $e_1,\ldots,e_g$ of pairwise disjoint
edges whose union is $X$ and 
such that $g$ is the height of $I(\mathcal{C})$.
\end{definition}

A set of pairwise disjoint edges is called {\it independent} and a set
of independent edges of $\mathcal{C}$ whose union is $X$ is called a
{\it perfect matching\/}. A clutter $\mathcal{C}$ satisfies the {\it
K\"onig property\/} if  
the maximum number of independent edges of $\mathcal{C}$ equals
the height of $I(\mathcal{C})$. It is rapidly seen that a 
clutter with a perfect matching of K\"onig type has 
the K\"onig property. In Lemma \ref{um+konig-kt} we show the converse
to be true for unmixed clutters. For uniform clutters, it is easy to
check that if $\mathcal{C}$ has the K\"onig property and a perfect
matching, then the perfect matching is of K\"onig type. However the
next example shows that this converse fails in general. 

\begin{example} Consider the clutter $\mathcal{C}$ with vertex set 
$X=\{x_1,\ldots,x_9\}$ whose edges are 
$$
\begin{array}{lllll}
e_1=\{x_1,x_2\},&e_2=\{x_3,x_4,x_5,x_6\},&e_3=\{x_7,x_8,x_9\},&  &
\\ 
f_4=\{x_1,x_3\},&f_5=\{x_2,x_4\},&f_6=\{x_5,x_7\},&f_7=\{x_6,x_8\}.&
\end{array}
$$
The edges $e_1,e_2,e_3$ form a perfect matching, 
$f_4,f_5,f_6,f_7$ are independent edges, and ${\rm ht}\, 
I(\mathcal{C})=4$. Thus $\mathcal{C}$ has the K\"onig property, but
$\mathcal{C}$ has no perfect matching of K\"onig type.  
\end{example}

\begin{lemma}\label{um+konig-kt} 
If $\mathcal{C}$ is an unmixed clutter with the 
K\"onig property and without isolated vertices, then $\mathcal{C}$ has
a perfect matching of K\"onig type.
\end{lemma}

\begin{proof} Let $X$ be the vertex set of $\mathcal{C}$. There are
$e_1,\ldots,e_g$ independent edges of $\mathcal{C}$, where $g$ is the
height of $I(\mathcal{C})$. If 
$e_1\cup\cdots\cup e_g\subsetneq X$, pick $x_r\in 
X\setminus(e_1\cup\cdots\cup e_g)$.
Since the vertex $x_r$ occurs in some edge of $\mathcal{C}$, there is
a minimal 
vertex cover $C$ containing $x_r$. Thus  
using that $e_1,\ldots,e_g$ are mutually disjoint we 
conclude that $C$ contains at least $g+1$ vertices, 
a contradiction. 
\end{proof}

\noindent {\it Notation.\/} As usual, we will 
use  $x^a$ as an abbreviation for $x_1^{a_1} \cdots x_n^{a_n}$, 
where $a=(a_1,\ldots,a_n)\in \mathbb{N}^n$. The {\it support\/} of a 
monomial $x^a=x_1^{a_1}\cdots x_n^{a_n}$ is given by ${\rm
supp}(x^a)= \{x_i\, |\, a_i>0\}$. 

\begin{proposition}\label{sas-phoenix07} Let $\mathcal{C}$ be an
unmixed clutter with a perfect 
matching $e_1,\ldots,e_g$ of K\"onig type 
and let $C_1,\ldots,C_r$ be any collection 
of minimal vertex covers of $\mathcal{C}$. If $\mathcal{C}'$ is 
the clutter associated to $I'=\cap_{i=1}^r(C_i)$, then $\mathcal{C}'$
has a perfect matching $e_1',\ldots,e_g'$ of K\"onig type such that:
{\rm(a)} 
$e_i'\subset e_i$ for all $i$,  
and {\rm (b)} every vertex of $e_i\setminus e_i'$ is isolated in 
$\mathcal{C}'$. 
\end{proposition}

\begin{proof} We denote the minimal set of generators 
of the ideal $I=I(\mathcal{C})$ by $G(I)$. There are monomials 
$x^{v_1},\ldots,x^{v_g}$ in $G(I)$ so that 
${\rm supp}(x^{v_i})=e_i$ for $i=1,\ldots,g$. Since $x^{v_i}$ is in 
$I$ and $I\subset I'$, there is $e_i'\subset e_i$ such 
that $e_i'$ is an edge of $\mathcal{C}'$. Let $x$ be any vertex in
$e_i\setminus e_i'$. If $x$ is not isolated in
$\mathcal{C}'$, there would a minimal vertex cover $C_k$ of
$\mathcal{C}'$ containing $x$. As $C_k$ contains a vertex of $e_j'$ 
for each $1\leq j\leq g$ and since $e_1',\ldots,e_g'$ are pairwise 
disjoint, we get that $C_k$ contains at least $g+1$
vertices, a contradiction. Thus (a) and (b) are satisfied. Clearly $g$
is the height of $I'$ by construction of $I'$. Let $X'$ be the vertex
set of $\mathcal{C}'$. To finish the proof we need only show that 
$X'=e_1'\cup\cdots\cup e_g'$. Let $x\in X'$, then 
$x\in e_i$ for some $i$ and $x$ belongs to at least one edge of
$\mathcal{C}'$. By part (b) we get that $x\in e_i'$, as
required. 
\end{proof}

\begin{remark}\label{may13-07} Let $C_1,\ldots,C_p$ be the minimal
vertex covers of 
$\mathcal{C}$. Since $I(\mathcal{C})$ is equal to $\cap_{i=1}^p(C_i)$, 
one has $(I(\mathcal{C})\colon x_j)=\cap_{x_j\notin C_i}(C_i)$ for any
vertex $x_j\notin I(\mathcal{C})$. Under the assumptions of
Proposition~\ref{sas-phoenix07} 
we get that $\mathcal{C}/x_j$ has a perfect matching
$e_1',\ldots,e_g'$ satisfying (a) and (b).
\end{remark}

\begin{lemma}\label{c-minor} Let $\mathcal{C}$ be an unmixed clutter
with a perfect  
matching $e_1,\ldots,e_g$ of K\"onig type and let $I=I(\mathcal{C})$.
If $e_1=\{x_1,\ldots,x_r\}$, then 
$$
\bigcap_{x_1 \in C_i}(C_i) = \left(((\cdots(((I\colon x_2)\colon
x_3)\colon 
x_4)\cdots)\colon x_{r-1})\colon x_r\right),
$$
where $C_1,\ldots,C_p$ are the minimal vertex covers 
of $\mathcal{C}$. 
\end{lemma}

\begin{proof} Let $I'$ denote the ideal on the right hand side of the
  equality. Then $I'$ is obtained from $I$ by 
making $x_i=1$ for $i=2,\ldots,r$, i.e., if $x^{v_1},\ldots,x^{v_q}$
generate $I$ and we make $x_i=1$ for $i=2,\ldots,r$ in
$x^{v_1},\ldots,x^{v_q}$, we obtain a generating set of $I'$. 
Notice that 
$I' = (I\colon x_2\cdots x_r)$ by the definition of the colon
operation. Take a 
monomial $x^a=x_1^{a_1}x_{r+1}^{a_{r+1}}\cdots x_n^{a_n}$ in $I'$.
We may assume $a_1=0$, otherwise $x^a$ is already in the left
hand side. Then $x_2\cdots x_rx_{r+1}^{a_{r+1}}\cdots x_n^{a_n}$ is in
$I$. Let $C_i$ be any minimal vertex cover of $ \mathcal{C}$
containing $x_1$. Observe that $C_i$ cannot contain $x_j$ 
for $2\leq j\leq r$. Indeed if $x_j\in C_i$ for some $2\leq j\leq r$,
then $C_i$ would contain $\{x_1,x_j\}$ plus at least one vertex of
each edge 
in the collection $e_2,\ldots,e_g$, a contradiction because $C_i$ has 
exactly $g$ vertices. Hence, using that 
$x_2\cdots x_r x_{r+1}^{a_{r+1}}\cdots x_n^{a_n}$ is in $I$, we
get that $x_{r+1}^{a_{r+1}}\cdots x_n^{a_n}$ is in $(C_i)$.
Consequently 
$x^a$ is in the left hand side of the equality. Conversely let $x^a$ be a 
minimal generator in the left hand side of the equality. Then $x^a \in
(C_i)$ whenever $x_1 \in C_i$. If $x_1 \not\in C_i$, then $x_2\cdots x_r
\in (C_i)$ since $C_i$ covers $e_1$. Thus $x^ax_2\cdots x_r \in (C_i)$ for
all $i$, and so $x^ax_2\cdots x_r \in \cap_{i=1}^p(C_i)=I$. Thus $x^a$
is in the right hand side of the equality.
\end{proof}

\begin{definition}
A simplicial complex $\Delta$ is {\it shellable\/}  if 
the facets (maximal faces) of $\Delta$ can be ordered
$F_1,\ldots,F_s$ such that  
for all $1\leq i<j\leq s$, there 
exists some $v\in F_j\setminus F_i$ and some 
$\ell\in \{1,\ldots,j-1\}$ with $F_j\setminus F_\ell= \{v\}$. We call
$F_1,\ldots,F_s$ a {\it shelling} of $\Delta$. 
\end{definition}

The above definition of shellable is due to
Bj\"orner and Wachs \cite{BW}.  Originally, the definition of
shellable also required that 
the simplicial complex be pure, that is, all the facets have same dimension.
We will say $\Delta$ is {\it pure shellable} if it also satisfies this
hypothesis. Because $I = I(\mathcal C)$ is a square-free monomial ideal,
it also corresponds to a simplicial complex via 
the Stanley-Reisner correspondence \cite{Stanley}.  
We let $\Delta_{\mathcal C}$ represent 
this simplicial complex. 
Note that $F$ is a facet of $\Delta_{\mathcal C}$ if and only if 
$X\setminus F$ is a minimal vertex cover of $\mathcal C$. 
For use below we say $x_i$ is a {\it free variable\/} (resp. {\it free
vertex}) of $I$ (resp. ${\mathcal C}$) if $x_i$
only appears in one of the monomials of $G(I)$ (resp. in
one of the edges of $\mathcal C$), where $G(I)$ denotes the 
minimal set of generators of the monomial ideal $I=I(\mathcal{C})$. 

If $\mathcal{C}$ has the {\it free vertex property}, i.e., 
all minors of $\mathcal{C}$ have a free vertex, then 
$\Delta_\mathcal{C}$ is shellable \cite{bipartite-scm}. 
We complement this result by showing that if all c-minors have a
free vertex and $\mathcal{C}$ is unmixed, then $\Delta_\mathcal{C}$ is
shellable. 

\begin{theorem}\label{c-minors} Let $\mathcal{C}$ be a clutter with a perfect
matching $e_1,\ldots,e_g$ of K\"onig type. If all c-minors of
$\mathcal C$ have a free  
vertex and $\mathcal C$ is unmixed, then $\Delta_{\mathcal C}$ is pure
shellable. 
\end{theorem}

\begin{proof} The proof is by induction on the number of vertices. We
may assume that $\mathcal{C}$ is a non-discrete clutter, i.e., it
contains an edge with at least two vertices. 
Let $z$ be a free vertex of $\mathcal{C}$ and let $C_1,\ldots,C_p$ be
the minimal vertex covers of $\mathcal{C}$. We may also assume that 
$z\in e_m$ for some $e_m=\{z_1,\ldots,z_r\}$, with $r\geq 2$. 
For simplicity of notation assume 
that $z=z_1$ and $m=g$. Consider 
the clutters $\mathcal{C}_1$ and $\mathcal{C}_2$ associated with
\begin{equation}\label{oct14-07}
I_1=\bigcap_{z_1\notin C_i}(C_i)\ \mbox{ and }\ I_2=\bigcap_{z_1\in
C_i}(C_i)
\end{equation}
respectively. By Proposition~\ref{sas-phoenix07}, the clutter
$\mathcal{C}_2$ has a perfect matching 
$e_1',\ldots,e_g'$ of K\"onig type such that: (a) $e_i'\subset e_i$
for all $i$,  
and (b) every vertex $x$ of $e_i\setminus e_i'$ is isolated in 
$\mathcal{C}_2$, i.e., $x$ does not occur in any edge of
$\mathcal{C}_2$. In particular all vertices of $e_g\setminus\{z_1\}$
are isolated vertices of 
$\mathcal{C}_2$. Similar statements hold for $\mathcal{C}_1$ because
of Proposition~\ref{sas-phoenix07}. By
Lemma~\ref{c-minor} and  
Remark~\ref{may13-07} we get
$$
I_1=(I\colon z_1)\ \mbox{ and }\ I_2=\left(((\cdots(((I\colon
z_2)\colon z_3)\colon 
z_4)\cdots)\colon z_{r-1})\colon z_r\right),
$$
that is, $\mathcal{C}_1=\mathcal{C}/z_1$ 
and 
$\mathcal{C}_2=\mathcal{C}/\{z_2,\ldots,z_r\}$. Hence the ideals   
$I_1$ and $I_2$ are c-minors of $I$. The number of vertices of
$\mathcal{C}_i$ is less than that of $\mathcal{C}$ for $i=1,2$.  
Thus $\Delta_{\mathcal{C}_1}$ and
$\Delta_{\mathcal{C}_2}$ are  
shellable by the induction hypothesis. Consider the clutter
$\mathcal{C}_i'$ whose edges are the edges of $\mathcal{C}_i$ and
whose vertex set is $X$. The minimal vertex covers of $\mathcal{C}_i'$
are exactly the minimal vertex covers of $\mathcal{C}_i$. Thus it
follows that $\Delta_{\mathcal{C}_i'}$ is shellable for $i=1,2$. Let
$F_1,\ldots,F_s$ be the facets 
of $\Delta_{\mathcal C}$ that contain $z_1$ and let $G_1,\ldots,G_t$
be the facets 
of $\Delta_{\mathcal C}$ that do not contain $z_1$. Notice that the
edge ideals 
of $\mathcal{C}_i$ and $\mathcal{C}_i'$ coincide, the vertex set 
of $\mathcal{C}_i'$ is equal to the vertex set of $\mathcal{C}$, 
and $I=I_1\cap I_2$. Hence from Eq.~(\ref{oct14-07}) we get that
$F_1,\ldots,F_s$ are the facets of $\Delta_{\mathcal{C}_1'}$ and
$G_1,\ldots,G_t$  are the facets of $\Delta_{\mathcal{C}_2'}$. By the
induction 
hypothesis we may assume $F_1,\ldots,F_s$ is a shelling of
$\Delta_{\mathcal{C}_1'}$ and $G_1,\ldots,G_t$ is a shelling of
$\Delta_{\mathcal{C}_2'}$. We now prove that
$$
F_1,\ldots,F_s,G_1,\ldots,G_t
$$
is a shelling of $\Delta_{\mathcal{C}}$. We need only show that 
given $G_j$ and $F_i$ there is $v\in G_j\setminus F_i$ and $F_\ell$
such that $G_j\setminus F_\ell=\{v\}$. We can write 
$$
G_j=X\setminus C_j\ \mbox{ and }\ F_i=X\setminus C_i,
$$
where $C_j$ (resp. $C_i$) is a minimal vertex cover of $\mathcal C$
containing $z_1$ (resp. not containing $z_1$). Notice that
$z_2,\ldots,z_r$ are not in $C_j$ because $e_1,\ldots,e_g$ is a
perfect matching and $|C_j|=g$. Thus $z_2,\ldots,z_r$ are in $G_j$.
Since $z_1\in F_i$ and $F_i$ cannot contain the edge $e_g$, there is a
$z_k$ so  
that $z_k\notin F_i$ and $k\neq 1$. Set $v=z_k$ and
$F_\ell=(G_j\setminus\{z_k\})\cup\{z_1\}$. Clearly $F_\ell$ is an 
independent vertex set because $z_1$ is a free vertex in $e_g$ and
$G_j$ is an 
independent vertex set. Thus $F_\ell$ is a facet because
$\mathcal{C}$ is 
unmixed. To complete the proof observe that $G_j\setminus
F_\ell=\{z_k\}$. 
\end{proof}

For use below we set $x_e=\prod_{x_i\in e}x_i$ for any $e\subset X$.
Next we give a characterization of the unmixed
property of $\mathcal{C}$. This characterization can be formulated 
combinatorially or algebraically.

\begin{theorem}\label{unmixed-clutter-pm} 
Let $\mathcal C$ be a clutter with a perfect matching $e_1,\ldots,e_g$
of K\"onig type and let $I=I(\mathcal{C})$ be its edge ideal. 
Then the following are equivalent:
\begin{itemize}
\item[(a)] $\mathcal C$ is unmixed.
\item[(b)] For any two edges $e \neq e'$ and for any two distinct vertices 
$x\in e$, $y\in e'$ contained in some $e_i$, one has that 
$(e\setminus\{x\})\cup (e'\setminus\{y\})$ contains an edge. 
\item[(c)] For any two edges $e\neq e'$ and for any $T\subset e_i$
such that $x_T$ divides $x_ex_{e'}$, one has that 
${\rm supp}(x_ex_{e'}/x_T)$ contains an edge.
\item[(d)] For any two edges $e\neq e'$ and for any $e_i$, 
$(x_ex_{e'}\colon x_{e_i}) \subset I$.
\item[(e)] $I=(I^2\colon x_{e_1})+\cdots+(I^2\colon x_{e_g})$.
\end{itemize}
\end{theorem}

\begin{proof} 
(a) $\Rightarrow$ (c): We may assume $i=1$. Let $T$ be a
subset of $e_1$ such that $x_T$ divides $x_ex_{e'}$. If $T \subset
e$, then 
$e'$ is an edge contained in $S={\rm supp}(x_e x_{e'}/x_T)$ and there
is nothing 
to show. The proof is similar
if $T \subset e'$. So we can define $T_1= e \cap T$ and $T_2 = T
\setminus T_1$ 
and we may assume neither $T_1$ nor $T_2$ is empty. Note that $T_1
\subset e$ and 
$T_2 \subset e'$. In fact, $T_2 \subset T \cap e'$, but equality does
not necessarily 
hold. Notice that $S =(e \setminus T_1) \cup (e' \setminus T_2)$. 
If $S$ does not contain an edge, its complement contains a minimal
vertex 
cover $C$. We use $c$ to denote complement. Then 
$$
C \subset X\setminus S = S^c =
(e \setminus T_1)^c \cap (e' \setminus T_2)^c = 
(e^c \cup T_1) \cap (e'^c \cup T_2).$$ 
Now $C \cap e \not= \emptyset$,  so there is an $x \in C \cap
e$. Then $x \in e^c \cup T_1$. This
forces $x \in T_1$. Similarly there is a $y \in C \cap e'$, and so 
$y \in e'^c \cap 
T_2$. Thus $y \in T_2$.  By the definition of $T_2$, $x \not= y$. To
derive a contradiction 
pick $z_k\in e_k\cap C$ for $k\geq 2$ and notice that
$x,y,z_2,\ldots,z_g$ is a set of $g+1$ distinct vertices in $C$, which
is impossible because $\mathcal{C}$ is unmixed.

(c) $\Rightarrow$ (b): Let $x \in e$ and $y \in e'$ be two distinct
vertices contained in some $e_i$. Let $T=\{x,y\}$. Then $x_T$ divides
$x_ex_{e'}$ and 
$$
S={\rm supp}(x_ex_{e'}/x_T) \subset
(e\setminus\{x\})\cup (e'\setminus\{y\}).
$$ 
By (c), $S$ contains an
edge. Thus $(e\setminus\{x\})\cup (e'\setminus\{y\})$ contains an edge.    

(b) $\Rightarrow$ (a):  Let $C$ be a minimal vertex cover of
$\mathcal{C}$. 
Since the matching 
is perfect, there is a partition:
$$
C=(C\cap e_1)\cup\cdots\cup (C\cap e_g).
$$
Hence it suffices to prove that $|C\cap e_i|=1$ for all $i$. We
proceed by contradiction. For simplicity of notation assume $i=1$ and
$|C\cap e_1|\geq 2$. Pick $x\neq y$ in $C\cap e_1$. Since $C$
is minimal, there are edges $e,e'$ such that 
\begin{equation}\label{lolita-taco}
e\cap (C\setminus\{x\})=\emptyset\mbox{ and } e'\cap
(C\setminus\{y\})=\emptyset.
\end{equation}
Clearly $x\in e$, $y\in e'$, and $e\neq e'$ because $y\notin e$. 
Then by hypothesis the set $S=(e\setminus\{x\})\cup
(e'\setminus\{y\})$ contains an 
edge $e''$. Take $z\in e''\cap C$, then $z\in e\setminus\{x\}$ or 
$z\in e'\setminus\{y\}$, which is impossible by
Eq.~(\ref{lolita-taco}).

(c) $\Rightarrow$ (d): Let $x^a \in (x_ex_{e'}\colon x_{e_i})$ be a
monomial generator of the colon ideal. Then $x^ax_{e_i} = m
x_ex_{e'}$ for some 
monomial $m$. Let $T \subset e_i$ be maximal such that $x_T$ divides
$x_ex_{e'}$. Then $x_{e_i \setminus T}$ divides $m$, and $x^a=
(m/x_{e_i \setminus T})(x_ex_{e'}/x_T)$. Since ${\rm
supp}(x_ex_{e'}/x_T)$ contains an edge, we have $x_ex_{e'}/x_T \in
I$. Thus $x^a \in I$ as desired.      

(d) $\Rightarrow$ (c): Suppose $T \subset e_i$ is such that $x_T$
divides $x_ex_{e'}$. Then 
$$(x_ex_{e'}/x_T)x_{e_i} =
x_ex_{e'}x_{e_i\setminus T},$$
and so $(x_ex_{e'}/x_T) \in (x_ex_{e'}
\colon x_{e_i})\subset I$. Thus $(x_ex_{e'}/x_T)$ is a multiple of a
monomial generator of $I$. Hence ${\rm supp}(x_ex_{e'}/x_T)$ contains
an edge.      

(e) $\Rightarrow$ (d): If equality in (e) holds, then 
$(I^2\colon x_{e_i})=I$ for all $i$. Hence from the inclusion 
$(I^2\colon x_{e_i})\subset I$ we rapidly obtain that condition (d)
holds.

(d) $\Rightarrow$ (e): It suffices to verify that 
$(I^2\colon x_{e_i})=I$ for all $i$. Since $I$ is clearly contained 
in $(I^2\colon x_{e_i})$, we need only show the inclusion 
$(I^2\colon x_{e_i})\subset I$. Take $x^a\in (I^2\colon x_{e_i})$,
then $x^ax_{e_i}=mx_ex_{e'}$ for some edges $e,e'$ of $\mathcal{C}$
and some monomial $m$. If $e\neq e'$, then by hypothesis $x^a\in
(x_ex_{e'}\colon 
x_{e_i})\subset I$, i.e., $x^a\in I$. If $e=e'$, then
$x^ax_{e_i}=mx_e^2$. Thus $x_e$ divides $x^a$ because $x_{e_i}$ 
is a square-free monomial, but this means that $x^a\in I$, as required. 
\end{proof}

\begin{definition}\label{bala-clu} 
Let $A$ be the incidence matrix of a clutter $\mathcal{C}$. A clutter
$\mathcal{C}$ has a {\it cycle\/} of length $r$ if there is a square
sub-matrix of $A$ of 
order $r\geq 3$ with exactly two $1$'s in 
each row and column.  A clutter without odd cycles is 
called {\it balanced} and an acyclic clutter is called {\it totally
balanced\/}.
\end{definition}

This definition of cycle is equivalent to the usual definition of 
cycle in the sense of hypergraph theory
\cite{berge,hhtz}. All minors of a balanced clutter have the 
K\"onig property \cite{Schr2}. If $G$ is a graph, then $G$ is balanced
if and only if $G$ is bipartite and $G$ is totally balanced if and
only if 
$G$ is a forest. 

The following result extends---to clutters with the K\"onig
property---an unmixedness criterion of \cite{unmixed} valid 
for bipartite graphs. As a byproduct we obtain a full description of
all unmixed balanced  
clutters. 
 
\begin{corollary}\label{char-bal-um} Let $\mathcal{C}$ be a clutter
with the K\"onig property. Then $\mathcal{C}$ is unmixed if and only
if there is a perfect 
matching $e_1,\ldots,e_g$ of K\"onig type such that for any two edges
$e\neq e'$ and for any two distinct vertices $x\in e$, $y\in e'$
contained in some $e_i$, one has that 
$(e\setminus\{x\})\cup (e'\setminus\{y\})$ contains an edge.
\end{corollary}

\begin{proof} $\Rightarrow$) Assume that $\mathcal{C}$ is unmixed. 
By Theorem~\ref{unmixed-clutter-pm} it suffices to observe that 
any unmixed clutter with the 
K\"onig property and without isolated vertices has a perfect
matching of K\"onig type, see Lemma~\ref{um+konig-kt}. 

$\Leftarrow$) This implication follows at once from
Theorem~\ref{unmixed-clutter-pm}. 
\end{proof}

The following ordering property was shown to hold for the clutter of 
facets of any unmixed simplicial tree \cite[Remark~7.2,
Corollary~7.8]{faridijct}.  
The next result is a wide 
generalization of this fact because unmixed simplicial trees are 
acyclic \cite{hhtz}, being balanced they have the K\"onig property 
\cite[Theorem~83.1]{Schr2}, and by Lemma~\ref{um+konig-kt} they have a
perfect matching. 

\begin{theorem}\label{jun1-07} Let $\mathcal{C}$ be a clutter 
with a perfect matching $e_1,\ldots,e_g$ of K\"onig type. If 
$\mathcal{C}$ has no cycles of length $3$ or $4$ and 
$\mathcal{C}$ is unmixed, then for any two edges $f_1,f_2$ of
$\mathcal{C}$ and  
for any $e_i$, one has that 
$f_1\cap e_i\subset f_2\cap e_i$ or $f_2\cap e_i\subset f_1\cap e_i$.
\end{theorem}

\begin{proof} For simplicity assume $i=1$. We proceed
by contradiction. Assume there are  
$x_1\in f_1\cap e_1\setminus f_2\cap e_1$ and $x_2\in f_2\cap
e_1\setminus f_1\cap e_1$. As $\mathcal{C}$ is unmixed, by
Theorem~\ref{unmixed-clutter-pm}(b) there is an
edge $e$ of $\mathcal{C}$ such that
$$
e\subset(f_1\setminus\{x_1\})\cup(f_2\setminus\{x_2\})=(f_1\cup
f_2)\setminus\{x_1,x_2\}.
$$
Since $e\not\subset e_1$, there is $x_3\in e\setminus e_1$. Then
either $x_3\in f_1$ or $x_3\in f_2$. Without loss of generality we may
assume $x_3\in f_1\setminus e_1$. For use below we denote the
incidence matrix of 
$\mathcal{C}$ by $A$.  

Case(I): $x_3\in f_2$. Then the matrix 
\[
\begin{array}{cccc}
& f_1&f_2&e_1\\
x_1& 1&0&1\\
x_2& 0&1&1\\
x_3& 1&1&0\\
\end{array}
\]
is a submatrix of $A$, a contradiction.

Case(II): $x_3\notin f_2$. Notice that $e\not\subset f_1$, otherwise
$e=f_1$ which is impossible because $x_1\in f_1\setminus e$. Thus
there is $x_4\in e\setminus f_1$ and $x_4\in (e\cap f_2)\setminus
f_1$. 

Subcase(II.a): $x_4\in e_1$. Then the matrix 
\[
\begin{array}{cccc}
& f_1&e&e_1\\
x_1& 1&0&1\\
x_3& 1&1&0\\
x_4& 0&1&1\\
\end{array}
\]
is a submatrix of $A$, a contradiction.

Subcase(II.b): $x_4\notin e_1$. Then the matrix 
\[
\begin{array}{ccccc}
& f_1&e&f_2&e_1\\
x_1& 1&0&0&1\\
x_2&0&0&1&1\\
x_3& 1&1&0&0\\
x_4&0&1&1&0
\end{array}
\]
is a submatrix of $A$, a contradiction.
\end{proof}

Conversely, the above ordering property implies unmixedness. Note
that the assumption on the incidence matrix is not needed for this
implication.  

\begin{theorem}\label{jun1-07-1} Let $\mathcal{C}$ be a clutter 
with a perfect matching $e_1,\ldots,e_g$ of K\"onig type. If for any
two edges  $f_1,f_2$ of $\mathcal{C}$ and 
for any $e_i$, one has that 
$f_1\cap e_i\subset f_2\cap e_i$ or 
$f_2\cap e_i\subset f_1\cap e_i$, then $\mathcal{C}$ is unmixed.  
\end{theorem}

\begin{proof} To show that $\mathcal{C}$ is unmixed it suffices to
verify condition (b) of Theorem~\ref{unmixed-clutter-pm}. 
Let $f_1\neq f_2$ be two edges
and let $x\in f_1$, $y\in f_2$ be two distinct vertices contained in
some 
$e_i$. For simplicity we assume $i=1$. Set $B=(f_1\setminus\{x\})\cup
(f_2\setminus\{y\})$. Then $f_2\cap e_1\subset
f_1\cap e_1$ or $f_1\cap e_1\subset
f_2\cap e_1$. In the first case we have that 
$f_2\subset B$. Indeed let $z\in f_2$. If $z\neq y$, 
then $z\in f_2\setminus\{y\}\subset B$, and if $z=y$, 
then  $z\in f_2\cap e_1\subset f_1\cap e_1$ and $z\neq x$, i.e., 
$z\in f_1\setminus\{x\}\subset B$. In the second case $f_1\subset B$.
\end{proof}

\begin{proposition}\label{structure-theorem}
Let $\mathcal{C}$ be an unmixed clutter without cycles of length 
$3$ or $4$. If $e_1,\ldots,e_g$ is a perfect matching 
of $\mathcal{C}$ of K\"onig type, then $e_i$ has a free vertex for all $i$.
\end{proposition}

\begin{proof} Fix an integer $i$ in $[1,g]$. We may assume that $e_i$
has at least 
one non-free vertex. Consider the set of edges:
$$
\mathcal{F}=\{f\in E({\mathcal C})\vert\, e_i\cap f\neq\emptyset;\, 
f\neq e_i\}.
$$
By Theorem~\ref{jun1-07}, the edges of $\mathcal{F}$ can be 
listed as $f_1,\ldots,f_r$ so that they satisfy the 
inclusions 
$$
f_1\cap e_i\subset f_2\cap e_i\subset\cdots\subset f_r\cap
e_i\subsetneq e_i.
$$
Thus any vertex of $e_i\setminus(f_r\cap e_i)$ is a free vertex of
$e_i$.
\end{proof}

\begin{theorem}\label{jun2-07} Let $\mathcal{C}$ be an unmixed
clutter with a perfect  
matching $e_1,\ldots,e_g$ of K\"onig type. If $\mathcal{C}$ has no
cycles of length $3$ or $4$, 
then $\Delta_{\mathcal C}$ is pure shellable. 
\end{theorem}

\begin{proof} All hypothesis are preserved under contractions, i.e.,
under c-minors. This follows from Remark~\ref{may13-07} and the fact
that the incidence matrix of a contraction of $\mathcal{C}$ is a 
submatrix of the incidence matrix of $\mathcal{C}$. Thus by
Proposition~\ref{structure-theorem} any 
c-minor has a free vertex and the result follows from
Theorem~\ref{c-minors}. 
\end{proof}

\begin{theorem}\label{ordering+matching-shellable} 
Let $\mathcal{C}$ be a clutter with a perfect 
matching $e_1,\ldots,e_g$ of K\"onig type. If 
for any two edges $f_1,f_2$ of $\mathcal{C}$ and 
for any edge $e_i$ of the perfect matching, one has that 
$f_1\cap e_i\subset f_2\cap e_i$ or $f_2\cap e_i\subset f_1\cap e_i$,
then $\Delta_{\mathcal C}$ is pure shellable. 
\end{theorem}

\begin{proof} Notice the following two assertions: 
(i) $\mathcal{C}$ is an unmixed clutter, which 
follows from Theorem~\ref{jun1-07-1},
and (ii) $e_i$ has a free vertex for all $i$, which follows from the 
proof of Proposition~\ref{structure-theorem}. Thus by
Theorem ~\ref{c-minors} we need only show that any c-minor has a free
vertex. By (ii) it suffices to show that our hypotheses are closed
under 
contractions. Let $x$ be a vertex of $\mathcal{C}$ and let
$\mathcal{C}'=\mathcal{C}/x$. By Remark~\ref{may13-07}, 
we get that $\mathcal{C}/x$ has a perfect matching
$e_1',\ldots,e_g'$ satisfying: (a) $e_i'\subset e_i$ for all $i$, 
and (b) every vertex of $e_i\setminus e_i'$ is isolated in
$\mathcal{C}'$. Let $e,e'$ be two edges of $\mathcal{C}'$ and let
$e_i'$ be an edge of the perfect matching of $\mathcal{C}'$. There are
edges $f,f'$ of $\mathcal{C}$ such that one of the following is
satisfied: $e=f$ and $e'=f'\setminus\{x\}$, $e=f\setminus\{x\}$ and
$e'=f'$, $e=f\setminus\{x\}$ and $e'=f'\setminus\{x\}$, $e=f$ and
$e'=f'$. We may assume $f\cap e_i\subset f'\cap e_i$. To finish the
proof we now show that $e\cap e_i'\subset e'\cap e_i'$. Take $z\in
e\cap e_i'$. Then $z\in f\cap e_i$ and consequently $z\in f'\cap
e_i'$. Since $x\notin e_i'$, one has $z\neq x$. It follow that $z\in
e'\cap e_i'$. 
\end{proof}

Let $G$ be a graph and let $V$ be its vertex set. For use below
consider the 
graph $G\cup W(V)$ obtained from $G$ by
adding new vertices $\{y_i\,\vert\,x_i\in V\}$ and new
edges $\{\{x_i,y_i\}\,\vert\,x_i\in V\}$. The edges $\{x_i,y_i\}$ are
called {\it whiskers}.   The notion of a whisker was introduced in
\cite[p.~392]{ITG}. 

\begin{corollary} If $G$ is a graph and $G'=G\cup W(V)$, then
$\Delta_{G'}$ is pure shellable.
\end{corollary}

\begin{proof} It follows at once from
Theorem~\ref{ordering+matching-shellable}. Indeed if
$V=\{x_1,\ldots,x_n\}$, then $\{x_1,y_1\},\ldots,\{x_n,y_n\}$ is a
perfect matching of $G'$ satisfying the ordering condition in 
Theorem~\ref{ordering+matching-shellable}. 
\end{proof}

Recall that a clutter $\mathcal{C}$ is called {\it totally
balanced\/} if $\mathcal{C}$ is 
acyclic and that a graph $G$ is totally balanced if and only if
$G$ is a forest. 
Faridi \cite{Faridi} introduced the notion of a leaf for a simplicial
complex $\Delta$.  Precisely, a facet $F$ of $\Delta$ is a {\it leaf}
if $F$ is 
the only facet of $\Delta$, or there exists a facet $G \neq F$ in
$\Delta$ such that 
$F \cap F' \subset F \cap G$ for all facets $F' \neq F$ in $\Delta$.
A simplicial complex $\Delta$ is a {\it simplicial forest} if every
nonempty 
subcollection, i.e., a subcomplex whose
facets are also facets of $\Delta$, of $\Delta$ contains a leaf. 
Recently Herzog, Hibi, Trung and Zheng 
\cite[Theorem~3.2]{hhtz} showed that $\mathcal{C}$ is the clutter of
the facets of a simplicial forest if and only if $\mathcal{C}$ is a
totally balanced clutter. Soleyman Jahan and X. Zheng 
\cite[Corollary~3.1]{soleyman-zheng} showed that $\mathcal{C}$ is a
totally balanced clutter if and only  
if $\mathcal{C}$ satisfies the free vertex property. Altogether one
has:

\begin{proposition}{\rm(\cite{hhtz}, \cite{soleyman-zheng})}
Let $\mathcal{C}$ be a clutter. Then the following
conditions are equivalent
\begin{itemize}
\item[(a)] $\mathcal{C}$ is the clutter of the facets of a simplicial
forest. 
\item[(b)] $\mathcal{C}$ has the vertex free property.
\item[(c)] $\mathcal{C}$ is totally balanced.
\end{itemize}
\end{proposition}

Thus some of the results in \cite{faridijct} can be examined using
the combinatorial 
structure of totally balanced clutters
\cite[Chapter~83, p.~1439--1451]{Schr2}. Since totally balanced
clutters are acyclic and satisfy the K\"onig property \cite{Schr2}, 
the next result generalizes the Cohen-Macaulay criterion for trees
given in \cite[Theorem~2.4]{Vi2} and is a far reaching
generalization of Faridi's characterization of unmixed 
simplicial  trees \cite[Remark~7.2, Corollary~7.8]{faridijct}.

\begin{corollary}\label{char-tbc} Let $\mathcal{C}$ be a clutter with
the K\"onig property and without cycles of length $3$ or $4$. 
Then any of the following conditions are
equivalent\/{\rm :} 
\begin{itemize}

\item[(a)] $\mathcal{C}$ is unmixed.

\item[(b)] There is a perfect matching $e_1,\ldots,e_g$, 
$g={\rm ht}\, I(\mathcal{C})$, 
such that $e_i$ has a free vertex for all $i$, and for any 
two edges $f_1,f_2$ of $\mathcal{C}$ and 
for any edge $e_i$ of the perfect matching, one has that 
$f_1\cap e_i\subset f_2\cap e_i$ or $f_2\cap e_i\subset f_1\cap e_i$.

\item[(c)] $R/I(\mathcal{C})$ is Cohen-Macaulay.  

\item[(d)] $\Delta_{\mathcal C}$ is a pure shellable simplicial complex.
\end{itemize}
\end{corollary}

\begin{proof} Using Lemma~\ref{um+konig-kt}, 
Theorems~\ref{jun1-07} and \ref{jun1-07-1}, and
Proposition~\ref{structure-theorem}  
it follows readily that conditions (a) and 
(b) are equivalent. Since (a) is equivalent to (b), from 
Theorem~\ref{jun2-07} 
we get that (b) implies (d). That (d) implies
(c) and (c) implies (a) are well known properties, see for instance
\cite{Stanley,monalg}. 
\end{proof}

Next we give some applications to totally balanced clutters. We begin
by recalling some notions. Let $A$ be the incidence matrix of a
clutter $\mathcal{C}$. The matrix $A$ is called {\it
perfect\/} if the polytope defined by the system $x\geq 0;\
xA\leq\mathbf{1}$ is integral, i.e., it has only integral vertices.
Here $\mathbf{1}$ denotes the vector with all 
its entries equal to $1$. A {\it clique\/} of a graph $G$ is a subset
of the set 
of vertices that induces 
a complete subgraph. We will also call a complete subgraph of $G$ 
a clique. The {\it vertex-clique matrix\/} of a graph $G$ is the
$\{0,1\}$-matrix 
whose rows are indexed by the vertices of $G$ and whose columns are
the incidence vectors of the maximal cliques of $G$. 
Let $G$ be a graph. A {\it colouring\/} of the vertices of $G$ 
is an assignment of colours to the vertices of $G$ in such a way that
adjacent vertices 
have distinct colours. The {\it chromatic number\/} of $G$ 
is the minimal number of colours in a colouring of $G$.
A graph is {\it perfect\/} if for every induced subgraph $H$, the
chromatic 
number of $H$ equals the size of the largest complete subgraph 
of $H$. A clutter is called {\it uniform\/} if all its edges have 
the same size.

\begin{corollary}\label{structureoftrees} Let $\mathcal{C}$ be an
unmixed 
totally balanced
clutter with vertex set $X$. If $\mathcal{C}$ has no isolated
vertices and $g$ is the  
height of $I(\mathcal{C})$, then 
\begin{itemize}
\item[(a)] \cite[Theorem~6.8]{faridijct} $\mathcal{C}$ has a perfect 
matching $e_1,\ldots,e_g$ of K\"onig type such that $e_i$ has a
free vertex for all $i$.    

\item[(b)] \cite[Corollary~5.4]{bipartite-scm} 
$\Delta_\mathcal{C}$ is a pure shellable simplicial complex.

\item[(c)] $\mathcal{C}$ is the clutter of maximal cliques of a perfect
graph $G$. 

\item[(d)] The set of non-free vertices of $e_i$ is contained in a
maximal clique of $G$. 

\item[(e)] \cite[Proposition~5.8]{reesclu} If $\mathcal{C}$ is
uniform, there is a partition 
$X^1,\ldots,X^d$ of $X$ such that any edge of
$\mathcal{C}$ intersects any $X^i$ in exactly one vertex. 
\end{itemize}
\end{corollary}

\begin{proof} (a) and (b) follow at once from Corollary~\ref{char-tbc}. 
(c) Let $A$ be the incidence matrix of $\mathcal{C}$. 
According to \cite{berge-balanced}, \cite[Corollary~83.1a(vii),
1441]{Schr2} $\mathcal{C}$ is balanced if and only if every
submatrix of $A$ is perfect. By \cite{chvatal} there is a perfect
graph $G$ such 
that $A$ is the vertex-clique matrix of $G$, i.e., $\mathcal{C}$ is
the clutter of maximal cliques of $G$. (d) 
Consider the set
$$
\mathcal{G}=\{e_i\cap e\vert\, e\in E(\mathcal{C});\, e\neq e_i\}.
$$
By Theorem~\ref{jun1-07}, the sets in $\mathcal{G}$ can be 
listed in increasing order 
$$
f_1\cap e_i\subset f_2\cap e_i\subset\cdots\subset f_r\cap
e_i\subsetneq e_i,
$$
for some edges $f_1,\ldots,f_r$. Thus $e_i\cap f_r$ is exactly the set
of non-free vertices of $e_i$, and $f_r$ is the required maximal
clique. 
\end{proof}

We have included part (d) as one of the properties of totally balanced
uniform clutters because it serves as an introduction to the notion
of admissible clutter to be defined in the next section. 

\section{Admissible clutters with a perfect matching}\label{clutters-pm}

Let $X^1,\ldots,X^d$ and $e_1,\ldots,e_g$ be two partitions of a
finite set $X$ such that $|e_i\cap X^j|\leq 1$ 
for all $i,j$.  The variables of the polynomial ring $K[X]$ are
linearly ordered by: $x\prec y$ iff ($x\in X^i$, $y\in X^j$, $i<j$) or
($x,y\in X^i$, $x\in e_k$, $y\in e_\ell$, $k<\ell$). 

Let $e$ be a subset of $X$ of size $k$ such that $|e\cap X^i|\leq 1$
for all $i$. There are unique
integers $1\leq i_1<\cdots<i_k\leq d$ and integers $j_1,\dots,j_k\in [1,g]$ 
such that
$$
\emptyset\neq e\cap X^{i_1}=\{x_1\},\ \emptyset\neq e\cap
X^{i_2}=\{x_2\},\ldots,\ \emptyset\neq e\cap X^{i_k}=\{x_k\}
$$
and $x_1\in e_{j_1},\ldots,x_k\in e_{j_k}$. We say that $e$ is 
{\it admissible\/} if $i_1=1,i_2=2,\ldots,i_k=k$ and 
$j_1\leq\cdots\leq j_k$. We can represent an admissible set 
$e=\{x_1,\ldots,x_k\}$ as 
$e=x^1_{j_1}\cdots x^{k}_{j_k}$, i.e., $x_i=x^{i}_{j_i}$ 
and $x^{i}_{j_i}\in X^i\cap e_{j_i}$ for all $i$. A monomial
$x^a$ is admissible if ${\rm supp}(x^a)$ is admissible. A clutter 
$\mathcal{C}$ is called {\it admissible}
if $e_1,\ldots,e_g$ are edges of $\mathcal C$, 
$e_i$ is
admissible for all $i$, and all other edges are admissible sets not
contained in any of the $e_i$'s. We can think of $X^1,\ldots,X^d$ as
color classes that color the edges. 

\begin{lemma}\label{admissible-ht} 
If $\mathcal{C}$ is an admissible clutter, then 
$e_1,\ldots,e_g$ is a perfect matching of K\"onig type. 
\end{lemma}

\begin{proof} It suffices to prove that $g={\rm ht}\, I(\mathcal{C})$.
Clearly ${\rm ht}\, I(\mathcal{C})\geq g$ because any minimal vertex cover
of $\mathcal{C}$ must contain at least one vertex of each $e_i$ and
the $e_i$'s form a partition of $X$. For each $1\leq i\leq g$ there 
is $y_i=x_i^1$ so that $e_i\cap X^1=\{y_i\}$. Since the $e_i$'s form a
partition we have the equality  
$$
(e_1\cap X^1)\cup\cdots \cup (e_g\cap X^1)=X^1.
$$
Thus $|X^1|=g$. To complete the proof notice that $X^1$ is a vertex
cover of $\mathcal{C}$ because all edges of $\mathcal{C}$ are
admissible. This shows ${\rm ht}\, I(\mathcal{C})\leq g$,
as required.
\end{proof}

Admissible clutters with two color
classes $X^1$, $X^2$ are special types of bipartite graphs. They will
be examined in Section~\ref{pm-bipartite-konig}.

\begin{example} Consider the following balanced admissible clutter
with color  classes  $X^1,X^2,X^3$ and edges
$e_1,e_2,e_3,f_1,f_2,f_3$.
$$
\begin{array}{ccc}
\begin{array}{ccccc}
\   &\ &X^1&X^2&X^3\\ 
e_1&=&x_1&y_1&\\
e_2&=&x_2&y_2&z_2\\
e_3&=&x_3&y_3&
\end{array}
&\ \  \ \ \ \ \ &
\begin{array}{ccccc}
\   &\ &X^1&X^2&X^3\\ 
f_1&=&x_1&y_2&z_2\\
f_2&=&x_1&y_3&\\
f_3&=&x_2&y_3&
\end{array}
\end{array}
$$
This clutter is Cohen-Macaulay, and $e_1,e_2,e_3$ is a
perfect matching of K\"onig type.
\end{example}

\begin{example}
The uniform admissible clutters with three color classes 
$$
X^1=\{x_1,\ldots,x_g\},\ 
X^2=\{y_1,\ldots,y_g\},\ 
X^3=\{z_1,\ldots,z_g\}
$$
are, up to permutation of variables, 
exactly the clutters with a perfect matching $e_i=\{x_i,y_i,z_i\}$
for $i=1,\ldots,g$ such that all edges of $\mathcal{C}$ have the 
form $\{x_i,y_j,z_k\}$, with $1\leq i\leq j\leq k\leq g$.
\end{example}

\begin{example} Consider the following admissible uniform clutter
with edges 
$e_1,e_2,f_1$, perfect matching $e_1,e_2$, and color classes
$X^1,X^2,X^3$:
$$
\begin{array}{ccccc}
\   &\ &X^1&X^2&X^3\\ 
e_1&=&x_1&y_1&z_1\\
e_2&=&x_2&y_2&z_2\\
f_1&=&x_1&y_1&z_2\\
\end{array}
$$
This clutter is Cohen-Macaulay.
\end{example}

An examination of the Cohen-Macaulay and unmixed criteria for
bipartite graphs (see Theorem~\ref{herzog-hibi} and
Corollary~\ref{unmixed-bip}) suggests the following conjecture.  

\begin{conjecture}\label{problem} If $\mathcal{C}$ is an admissible
clutter and $\mathcal{C}$ is unmixed, then $I(\mathcal{C})$ is
Cohen-Macaulay.  
\end{conjecture}

This conjecture is true for admissible clutters with two color 
classes $X^1$, $X^2$ (see Theorem~\ref{herzog-hibi}) 
and has been verified in a large number of examples. 

Let $e_1,\ldots,e_g$ and $X^1,\ldots,X^d$ be as in the beginning 
of Section~\ref{clutters-pm}. 
Suppose $e_1,\ldots,e_g$ are admissible subsets of $X$. 
The clutter $\mathcal{C}$ on $X$ whose set of edges is:
$$
E({\mathcal C})=\left\{e\subset X\left\vert\begin{array}{l}
e_i\not\subset e\mbox{ for
}i=1,\ldots,g,\, e\mbox{ is admissible},\\ e\not\subset e'\mbox{ for
any admissible set } e'\neq e\end{array}
\right\}\right.\cup\{e_1,\ldots,e_g\}
$$
is called a {\it complete admissible clutter\/}. This clutter consists
of the maximal admissible sets with respect to inclusion. By 
Lemma~\ref{admissible-ht} we get that $e_1,\ldots,e_g$ is a perfect
matching of K\"onig type.

\begin{proposition}\label{complete-case-enrique} 
If $\mathcal{C}$ is a complete admissible 
clutter, then $\mathcal{C}$ is unmixed.
\end{proposition}

\begin{proof} To show that $\mathcal{C}$ is unmixed it suffices to
verify condition (b) of Theorem~\ref{unmixed-clutter-pm}. 
Let $e\neq e'$ be two edges of
$\mathcal{C}$ and let 
$x\neq y$ be two vertices such that $\{x,y\}\subset e_i$ for some $e_i$, 
$x\in e$, and $y\in e'$. Since $e,e',e_i$ are admissible we can write
$$
e=\{x_1,\ldots,x_k\},\ e'=\{y_1,\ldots,y_{k'}\},\
e_i=\{z_1,\ldots,z_{k''}\},
$$ 
where $x_i\in X^i$, $y_i\in X^i$, $z_i\in X^i$. Then there are
$i_1,i_2$ such that $x=x_{i_1}$, $y=y_{i_2}$, $x=z_{i_1}$, and
$y=z_{i_2}$. Without loss of generality we may assume  $i_1<i_2$. One
has $i_1<k$, because if $k=i_1$, then $e\subsetneq e\cup
\{z_{i_1+1},\ldots,z_{i_2}\}$ and the right hand side is admissible, a
contradiction. Set $f=\{y_1,\ldots,y_{i_1},x_{i_1+1},\ldots,x_k\}$.
Then  
$$
f\subset
e\setminus\{x\}\cup e'\setminus\{y\}.
$$
Thus to finish the proof we need only show that $f$ is an edge of 
$\mathcal{C}$. Since $y_{i_2}\in e_i$ and $x_{i_1}\in e_i$, then
$y_{i_1}\in e_\ell$ 
for some $\ell\leq i$ and $x_{i_1+1}\in
e_t$ for some $i\leq t$. Hence $f$ is admissible. 
Next we show that $f$ is maximal. Assume that $f$ is not maximal. 
Then there exists an admissible subset $f'$ that properly contains
$f$. Then there is $z\in 
f'\cap X^{k+1}$ and since $f\cup \{z\}\subset f'$, we get that 
$e\cup \{z\}=\{x_1,\ldots,x_k,z\}$ is admissible, 
but $e\subsetneq e\cup\{z\}$, a contradiction. Hence $f$ is maximal. 
\end{proof}

Suppose ${\mathcal C}$ is a clutter on the vertex set $X$ with a
perfect 
matching $e_1, \ldots , e_g$ where $g$ is the height of 
$I({\mathcal C})$, 
and let $X^1, \ldots , X^d$ be a 
partition of $X$ such that every edge of $\mathcal{C}$ intersects
each $X^i$ 
exactly once. If every maximal admissible subset of $X$ is an edge of
${\mathcal C}$ and these are the only edges of ${\mathcal C}$, then
we call 
${\mathcal C}$ a {\it complete admissible 
uniform clutter\/}. Note that a complete admissible uniform clutter
is in fact uniform with every edge having $d$ vertices. Also,
Proposition~\ref{complete-case-enrique}
holds and ${\mathcal C}$ is unmixed. 

\begin{theorem}\label{susan-shellable}
If ${\mathcal C}$ is a complete admissible uniform clutter, then the
simplicial complex generated by the edges of $\mathcal{C}$ is pure
shellable.  

\end{theorem}

\begin{proof} Order the variables of $K[X]$ as in the beginning of 
Section~\ref{clutters-pm}.
Since every monomial intersects each $X^i$ exactly
once, we can represent the edges of ${\mathcal C}$ as
$F_i=x^1_{i_1}x^2_{i_2}\cdots x^d_{i_d}$ where $x_{i_j}^i \in X^i\cap
e_{i_j}$ (example: $x^3_2 \in X^3 \cap
e_2$). Since $X^i \cap e_j$ has precisely one element for each $i,j$,
this notation is well-defined. Then we order the edges of
$\mathcal{C}$ lexicographically, that is $F_i=x^1_{i_1}x^2_{i_2}\cdots
x^d_{i_d} < F_j=x^1_{j_1}x^2_{j_2}\cdots 
x^d_{j_d}$ if the first nonzero entry of $(j_1, j_2, \ldots , j_d) -
(i_1, i_2, \ldots , i_d)={\bf j}-{\bf i}$ is positive. 
Under this
order, we show that 
${\mathcal C}$ is shellable. 

Suppose $F_i$ and $F_j$ are two edges of
${\mathcal C}$ with $F_i < F_j$. Suppose the first non-zero entry 
of ${\mathbf j} - {\mathbf i}$ is $j_t-i_t$. Then $1 \leq i_t < j_t$. 
Let $F_k = F_j \setminus \{ x^t_{j_t} \} \cup \{ x^t_{i_t} \}$ and 
let $v=x^t_{j_t}$. Since
$j_1=i_1\leq \cdots \leq j_{t-1}=i_{t-1} \leq i_t <j_t \leq j_{t+1}
\leq \cdots \leq j_d$ then $F_k$ is  maximal admissible,
$v\in F_j\setminus F_i$,
$F_k < F_j$ and  
$F_j\backslash F_k=\{v\}$ as required.
\end{proof}

\begin{example} The complete admissible uniform  clutter 
with three color classes 
$$
X^1=\{x_1,\ldots,x_g\},\ 
X^2=\{y_1,\ldots,y_g\},\ 
X^3=\{z_1,\ldots,z_g\}
$$
is the clutter $\mathcal{C}$ whose edge set is 
$E(\mathcal{C})=\{\{x_i,y_j,z_k\}\vert\, 1\leq i\leq
j\leq k\leq g\}$. Note that
$e_1=\{x_1,y_1,z_1\},\ldots,e_g=\{x_g,y_g,z_g\}$ 
is the perfect matching of $\mathcal{C}$.
\end{example}

The next example illustrates the construction of the 
lexicographical shelling used in the proof of
Theorem~\ref{susan-shellable}.

\begin{example} Let $\mathcal{C}$ be the complete admissible uniform  
clutter with color classes 
$X^1=\{x_1,x_2,x_3\}$, $X^2=\{y_1,y_2,y_3\}$,
$X^3=\{z_1,z_2,z_3\}$. Then the shelling of the simplicial complex 
generated by the  
edges of $\mathcal C$ is:
$$
\begin{array}{lllll}
F_1\ =\{x_1,y_1,z_1\}& < &F_2=\{x_1,y_1,z_2\} &<& 
F_3=\{x_1,y_1,z_3\}\ <\\
F_4\ =\{x_1,y_2,z_2\} &<& F_5 =\{x_1,y_2,z_3\} 
  &<& F_6=\{x_1,y_3,z_3\}\ < \\ 
F_7\ =\{x_2,y_2,z_2\} &<& F_8=\{x_2,y_2,z_3\} &<& 
F_9=\{x_2,y_3,z_3\}\ < \\ F_{10}=\{x_3,y_3,z_3\}. & & & &
\end{array}
$$
\end{example}

Let $\mathcal{C}$ be a clutter. The {\it Alexander dual\/} of
$\mathcal{C}$, denoted by $\Upsilon(\mathcal{C})$, is the clutter 
whose edges are the minimal vertex covers of $\mathcal{C}$. The edge
ideal of $\Upsilon(\mathcal{C})$ is called the Alexander dual of
$I(\mathcal{C})$. In combinatorial optimization the Alexander dual of
a clutter is referred to as the {\it blocker\/} of the clutter
\cite{Schr2}. 

\begin{lemma}\label{susan-duality}
If ${\mathcal C}$ is a complete admissible uniform clutter, then the
Alexander dual $\Upsilon(\mathcal{C})$ of ${\mathcal C}$ is also a
complete admissible uniform clutter. 
\end{lemma}

\begin{proof} Since ${\mathcal C}$ is
unmixed with covering number $g={\rm ht}\, I(\mathcal{C})$, then the 
Alexander dual  
is uniform with edges 
of size $g$. Note that $e_1, \ldots , e_g$ form a partition of the
vertices 
of the Alexander dual. Every minimal vertex cover of ${\mathcal C}$
must by definition intersect each $e_i$ at least once, and since
${\mathcal C}$ is unmixed
all minimal vertex covers have exactly $g$ elements, thus every edge
of $\Upsilon(\mathcal{C})$ intersects each 
$e_i$ exactly once. Also, $X^1, \ldots, X^d$ is a perfect matching of
$\Upsilon(\mathcal{C})$ since the $X^i$ partition the vertices and
since each edge of 
${\mathcal C}$ intersects each $X^i$ exactly once, $X^i$ is a minimal
vertex cover of ${\mathcal C}$, and thus an edge of the Alexander dual. 

Now since every minimal vertex cover of ${\mathcal C}$ has $g$ elements
and intersects $e_i$ exactly once for each $i$, all edges of the
Alexander dual have the form $M=x^{i_1}_1  x^{i_2}_2 \cdots
x^{i_g}_g$ where $1 
\leq i_t \leq d$ for all $1 \leq t \leq g$. To show
that the edges of $\Upsilon(\mathcal{C})$ are precisely the maximal
admissible 
subsets (with the $e_i$'s being the partition and the $X^i$'s the
perfect matching), 
we must show that $M$ is an edge of $\Upsilon(\mathcal{C})$ 
if and only if $i_1 \leq i_2 \leq
\cdots \leq i_g$. 

Suppose $M$ is as above and 
$1 \leq i_1\leq \cdots \leq i_g \leq d$. 
Suppose $F_j = x^1_{j_1}\cdots x^d_{j_d}$ is an edge of 
${\mathcal C}$. 
Then $F_j$ is admissible, so 
$1 \leq j_1 \leq \cdots \leq j_d \leq g$. 
We must show $M \cap F_j \not= \emptyset$. 
If $x^1_{j_1} \in M$ the intersection is not empty. 
Else, since $j_1 \in \{ 1 ,\ldots  g \},$ 
then $x^{i_{j_1}}_{j_1} \in M$ for some $i_{j_1} > 1$. 
Thus $i_t \geq 2$ for $t \geq j_1.$ Consider $x^2_{j_2}$. 
If $x^2_{j_2} \in M$, done. Else $i_t \geq 3$ for $t \geq j_2$. 
Since $i_d \leq d$, this process must stop with an element 
in the  intersection of $M$ and $F_j$, or $i_t = d$ for all 
$t \geq j_s$ for some $s.$ If $i_t =d$ for $t \geq j_s$, 
then since $j_s \leq g$ and $j_d \geq j_s$, we have 
$x^d_{j_d} \in F_j \cap M$
and thus the intersection is not empty and so $M$ is 
a minimal vertex cover of ${\mathcal C}$ and so an  edge of the
Alexander dual. 

Now suppose $M$ is as above, but $i_t > i_s$ for some $t<s.$ Choose
$t$ and $s$ so that $i_j < i_t$ for $j<t$ and $i_{\ell} \geq i_t$ for
$t<\ell < s.$ Define $F= x^1_t \cdots x^{i_{t}-1}_t x^{i_t}_s \cdots
x^d_s$. Then since $t<s,$ $F$ is maximal admissible and so an edge of
${\mathcal C}$. But $M\cap e_t = \{x^{i_t}_t\}$ and
$M\cap e_s = \{x^{i_s}_s\}$. 
Now $x^{i_t}_t \not\in F \cap e_t = \{x^1_t, \ldots x^{i_{t}-1}_t \}$
and since $i_s < i_t,$ $x^{i_s}_s \not\in F \cap e_s.$ Thus 
$F \cap M = \emptyset.$ Thus $M$ is not a vertex cover of ${\mathcal
C}$ and so 
is not an edge of the Alexander dual.
\end{proof}

\begin{lemma}\label{complete-complement}
If ${\mathcal C}$ is a complete admissible uniform clutter, 
then the simplicial complex $\Delta_{\Upsilon(\mathcal{C})}$ 
generated by $\{X\setminus F\vert\, F\in E(\mathcal{C})\}$ is pure 
shellable. 
\end{lemma}

\begin{proof} Let $F_1, \ldots F_r$ be the shelling of the edges of 
${\mathcal C}$ defined in Theorem \ref{susan-shellable}. Let 
$G_1 = X \setminus F_1,\ldots,G_r = X \setminus F_r$ 
be the facets of $\Delta_{\Upsilon(\mathcal{C})}$. We claim that 
$G_1,\ldots,G_r$ is the desired shelling. 
Suppose $G_i < G_j$. Then $F_i <
F_j$. Using the notation defined in Theorem \ref{susan-shellable},
let $v=x^t_{j_t}$ and define $u=x^t_{i_t}$. 
Then $u \in G_j\setminus G_i$ and $G_j\setminus G_k = \{ u \}$ as
required. 
\end{proof}

\begin{theorem}\label{complete-CM}
If ${\mathcal C}$ is a complete admissible uniform clutter, then
$R/I({\mathcal C})$ is a Cohen-Macaulay ring with a $d$-linear
resolution and $|E(\mathcal{C})|=\binom{d+g-1}{g-1}$.
\end{theorem}

\begin{proof} Consider the clutter $\Upsilon(\mathcal{C})$ 
of minimal vertex covers of $\mathcal{C}$. 
By Lemma \ref{complete-complement} and Lemma \ref{susan-duality} we
have that  $\Delta_{\Upsilon(\mathcal{C})}$ is pure shellable. 
Now recall that the Stanley-Reisner ideal of
$\Delta_{\Upsilon(\mathcal{C})}$ is $I(\Upsilon(\mathcal{C}))$ and
that  $I(\Upsilon(\mathcal{C}))$ is the Alexander dual of 
$I(\mathcal{C})$. Thus $R/I(\Upsilon(\mathcal{C}))$ is 
Cohen-Macaulay, 
and by \cite{ER} the ideal $I(\mathcal{C})$ has a linear resolution. 
Since the Alexander dual of a complete admissible uniform
clutter is also a complete admissible uniform clutter and since 
$\Upsilon(\Upsilon(\mathcal{C}))=\mathcal{C}$ it follows  
that $R/I(\mathcal{C})$ is Cohen-Macaulay. The formula for the 
number of edges of $\mathcal{C}$ follows from the 
explicit formula given in \cite{HKu} for the Betti numbers of a
Cohen-Macaulay ideal with a linear resolution. 
\end{proof}

Let $\mathcal{C}$ be a complete admissible uniform clutter. For each
edge  
$e=x^1_{j_1}x^2_{j_2}\cdots x^d_{j_d}$ of $\mathcal{C}$ consider all 
pairs $(x_{j_i}^i,x_{j_k}^k)$ with $i<k$ and consider the union of all
these pairs with $e=x^1_{j_1}x^2_{j_2}\cdots x^d_{j_d}$ running
through all edges of $\mathcal{C}$. This defines 
a poset $(P,\prec)$ on $X$ whose comparability graph $G$ is defined 
by all the unordered pairs $\{x_{j_i}^i,x_{j_k}^k\}$. The graph $G$
is perfect \cite[Corollary~{66.2a}]{Schr2} and any d-minor of the
clutter of maximal cliques of $G$ satisfies the K\"onig property. 
This follows from a variant of Dilworth's decomposition theorem
\cite[Theorem~14.18]{Schr2}. In the terminology of \cite{bonomo} $G$ 
is clique-perfect.

\begin{corollary} If $G'$ is the complement of the comparability
graph $G$ defined above, then $R/I(G')$ is Cohen-Macaulay.
\end{corollary}

\begin{proof} Notice that $\Delta_{G'}
=\{\mathcal{K}_r\vert\,
\mathcal{K}_r\mbox{ is a clique of }G\}=\mathcal{O}(P)$, where
$\mathcal{O}(P)$ 
is the order complex of $P$. Since the maximal faces 
of $\mathcal{O}(P)$ are precisely the edges
of $\mathcal{C}$, by Theorem~\ref{susan-shellable}, 
we obtain that $\mathcal{O}(P)$ is a pure shellable
complex whose Stanley-Reisner ring is $R/I(G')$. Hence $R/I(G')$ 
is Cohen-Macaulay. 
\end{proof}

Let $\mathcal{C}$ be a clutter and let $x^{v_1},\ldots,x^{v_q}$ be 
the minimal set of generators of $I(\mathcal{C})$. Consider the ideal
$I^*=(x^{w_1},\ldots,x^{w_q})$, where $v_i+w_i=(1,\ldots,1)$.
Following the terminology of matroid
theory we call $I^*$ the {\it dual\/} of $I$. Recall that $I^*$ has
 {\it linear quotients\/} if there is an ordering of the generators
 $x^{w_1},\ldots,x^{w_q}$ such that 
$$
((x^{w_1},\ldots,x^{w_{i-1}})\colon (x^{w_i})) = 
(x_{i_1},\ldots,x_{i_t})
$$
for $i=2,\ldots,q$, i.e., all colon ideals are generated by subsets
of the 
set of variables. If $I^*$ has linear quotients and all $x^{w_i}$
have the same degree, then $I^*$ has a linear resolution 
(see \cite[Lemma~5.2]{Faridi}, \cite{zheng-ca}).

\begin{corollary} If $\mathcal{C}$ is a complete admissible 
uniform clutter, then $I(\mathcal{C})^*$ has linear quotients.
\end{corollary}

\begin{proof} Let $x^{v_1},\ldots,x^{v_q}$ be the minimal set of
generators of $I=I(\mathcal{C})$ and let $F_i={\rm supp}(x^{v_i})$ for
$i=1,\ldots,q$. By Theorem~\ref{susan-shellable}, we may assume that 
$F_1,\ldots,F_q$ is a shelling for the 
simplicial complex $\langle F_1,\ldots,F_q\rangle$ generated by the 
$F_i$'s. Thus according to \cite[Theorem
1.4(c)]{hhz-ejc} the ideal $I^*=(x_{F_1^c},\ldots,x_{F_q^c})$ 
has linear quotients, where $F_k^c=X\setminus F_k$ and 
$\textstyle x_{F_k^c}=\prod_{x_i\in F_k^c}x_i$.
\end{proof}

We may also redefine the notion of admissible monomial to 
allow ``gaps''. This can be done as follows. 
Let $S=\{ x_1, \ldots , x_s \} $ be a subset of $X$ of size $s$ such
that $|S\cap X^i|\leq 1$ for all $i$. 
There are $k_1,\ldots,k_s$ and $j_1,\ldots,j_r$ such that 
$x_\ell\in X^{k_\ell}$ and $x_i\in e_{j_i}$ for all $i,\ell$. The set
$S$ is called {\it admissible} if 
$j_1\leq \cdots\leq j_r\leq g$ and  $k_1<\cdots< k_s$. A monomial
$x^a$ is admissible if ${\rm supp}(x^a)$ is admissible.

\begin{example} Consider the following clutter with edges
$e_1,e_2,f_1,f_2$ and color classes
$X^1,X^2,X^3$
$$
\begin{array}{ccc}
\begin{array}{ccccc}
\   &\ &X^1&X^2&X^3\\ 
e_1&=&x_1&y_1&\\
e_2&=&&y_2&z_2
\end{array}
&\ \ \ \  \ \ &
\begin{array}{ccccc}
\   &\ &X^1&X^2&X^3\\ 
f_1&=& &y_1&z_2\\
f_2&=&x_1&y_2& 
\end{array}
\end{array}
$$
This clutter is unmixed, non-Cohen-Macaulay, has a perfect
matching $e_1,e_2$ of K\"onig type, 
and the height of $I(\mathcal{C})$ is two. Thus this example
shows that allowing gaps gives a negative answer to
Conjecture~\ref{problem}.
\end{example}

\section{Cohen Macaulay bipartite graphs and
shellability}\label{pm-bipartite-konig} 

Throughout this section we assume that $G$ is a bipartite graph with 
bipartition $V_1=\{x_1,\ldots,x_g\}$ and $V_2=\{y_1,\ldots,y_{g_1}\}$
and without isolated vertices. 

The following nice criterion of Herzog and Hibi classifies all
Cohen-Macaulay bipartite graphs.

\begin{theorem}[\rm\cite{herzog-hibi-crit}]\label{herzog-hibi} $G$ is 
Cohen-Macaulay bipartite graph 
if and only if $g=|V_1|=|V_2|$ and we can order the vertices such
that:  
{\rm ($\mathrm{h}_0$)} $\{x_i,y_i\}\in E(G)$ for $i=1,\ldots,g$, 
{\rm ($\mathrm{h}_1$)} if $\{x_i,y_j\}\in E(G)$, then $i\leq j$, and 
{\rm ($\mathrm{h}_2$)} if $\{x_i,y_j\}$ and $\{x_j,y_k\}$ are in
$E(G)$ and $i<j<k$,  
then $\{x_i,y_k\}\in E(G)$.
\end{theorem} 

The results of this section are inspired by this criterion. Below we 
study condition {\rm ($\mathrm{h}_1$)} and a variation of 
condition {\rm ($\mathrm{h}_2$)}. Observe that the uniform admissible
clutters with two color classes $X^1$,
$X^2$ (see Section~\ref{clutters-pm}) are exactly the bipartite
graphs that satisfy  {\rm 
($\mathrm{h}_0$)} and ($\mathrm{h}_1$).

Next we give a combinatorial characterization--suggested by condition
{\rm ($\mathrm{h}_2$)}-- of all unmixed bipartite graphs. 

\begin{corollary}[\rm\cite{unmixed}]\label{unmixed-bip} Let $G$ be a
bipartite graph. Then $G$ is unmixed 
if and only  there is a perfect matching $e_1,\ldots,e_g$ such that
  for any two edges $e\neq e'$ and for any two distinct vertices 
$x\in e$, $y\in e'$ contained in some $e_i$, one has that 
$(e\setminus\{x\})\cup (e'\setminus\{y\})$ is an edge.
\end{corollary} 

\begin{proof} It follows at once from Corollary~\ref{char-bal-um}
because bipartite graphs satisfy the K\"onig 
property \cite{Schr2}.   
\end{proof}

This corollary shows that condition {\rm ($\mathrm{h}_2$)} is in
essence 
an expression for the unmixed property of $G$, i.e., in
Theorem~\ref{herzog-hibi} we may assume that $G$ is unmixed instead
of  
assuming condition {\rm ($\mathrm{h}_2$)}. 

Let $\Delta_G$ be the Stanley-Reisner complex of $I(G)$. Its 
facets are the maximal independent (stable) sets of vertices of $G$.
Following 
\cite{duval} we define the $k${\it th} {\it pure skeleton\/} of 
$\Delta_G$ as:
\begin{eqnarray*}
\Delta_G^{[k]}&=&\langle\{F\in\Delta_G\vert\, k=|F|\}\rangle;\ \ 0\leq
k\leq\dim(\Delta_G)+1,
\end{eqnarray*}
where $\langle{\mathcal F}\rangle$ denotes the subcomplex generated by 
$\mathcal F$. Note that this simplicial complex is 
always pure. By an interesting result of Duval
\cite[Theorem~3.3]{duval} a simplicial complex $\Delta$ is 
sequentially Cohen-Macaulay if and only if $\Delta^{[k]}$ 
is Cohen-Macaulay for  $0\leq k\leq \dim(\Delta)+1$. In particular
$R/I(G)$ is Cohen-Macaulay if and only if $R/I(G)$ is sequentially
Cohen-Macaulay and $G$ is unmixed. Here we shall be
interested only in the pure skeleton of $\Delta_G$ of maximum
dimension.

The following result characterizes all bipartite 
graphs with a perfect matching that satisfy condition
($\mathrm{h}_1$). It gives a combinatorial 
description of the admissible uniform clutters with two color classes.

We come to the main result of this section.

\begin{theorem}\label{ordering-char} If $G$ is a bipartite graph with
a perfect matching
$e_1,\ldots,e_g$ such that $e_i=\{x_i,y_i\}$ for all $i$, 
then $\Gamma=\Delta_G^{[g]}$ is pure shellable if and only if we
can order $e_1,\ldots,e_g$ such that $\{x_i,y_j\}\in E(G)$ implies 
$i\leq j$.
\end{theorem}

\begin{proof} $\Leftarrow$) It suffices to show that
$\Gamma=\Delta_G^{[g]}$ 
is shellable because this simplicial complex is always pure. We
proceed by 
induction on $g$. 
Each facet of
$\Gamma$ contains exactly one vertex of each edge of the perfect 
matching. We set 
$$
A=\{y_i\vert\, x_i\in N(y_g)\};\ \ B=A\cup N(y_g)=\bigcup_{x_i\in
N(y_g)}\{x_i,y_i\},
$$
where $N(y_g)$ is the set of vertices of $G$ that are adjacent to
$y_g$. Consider the graph $G'=G\setminus B$, obtained from $G$ by
removing 
all vertices of $B$ and all edges incident with some vertex of $B$. 

Let $F_1'=\emptyset$ if $|A|=g$, in which case $G'=\emptyset$. 
Else let $F_1',\ldots,F_r'$ be the facets of 
$\Gamma'=\Delta_{G'}^{[\ell]}$
that do not 
intersect $N(A)$, where $\ell=g-|A|$. Here $N(A)$ denotes the 
{\it neighbor set} of $A$, i.e., the set of
vertices of $G$ that are adjacent to some vertex of $A$. 
We claim that 
$F_1=F_1'\cup A,\ldots,F_r=F_r'\cup A$ is the set of facets of
$\Gamma$ that contain $y_g$. 
First we show that $F_k$ is a
facet of $\Gamma$ for all $k$. If $F_k$ contains an edge
$e=\{x_i,y_j\}$, then  
$y_j\in A$ and $x_i\in F_k'$ because $A$ and $F_k'$ are independent.
Then $x_i\in N(A)$, a contradiction because $N(A)\cap F_k'=\emptyset$.
Hence $F_k$ is independent and it is a facet of $\Gamma$ because
$|F_k|=g$. Conversely, let $F$ be a facet of $\Gamma$ containing
$y_g$. Then $F\cap N(y_g)=\emptyset$, $A\subset F$, and $F\cap
N(A)=\emptyset$. Thus we can write $F=F'\cup A$, where $F'=F\setminus
A$ is a facet of 
$\Gamma'$ with $F'\cap N(A)=\emptyset$, as
required. By the induction hypothesis $\Gamma'$ is shellable. Next we
prove that $F_1',\ldots,F_r'$ is a shelling with the linear order induced
by the shelling of $\Gamma'$. Assume $F_i'<F_j'$. Since $\Gamma'$ is
shellable, there are $v\in F_j'\setminus F_i'$ and a facet $F'$ of
$\Gamma'$ such that $F'<F_j'$ and $F_j'\setminus F'=\{v\}$. It
suffices to prove that $F'$ does not intersect $N(A)$. If $F'\cap
N(A)\neq\emptyset$, pick $x_p$ in the intersection.  Then $x_p\notin
F_i'\cup F_j'$ because $F_i'$ and $F_j'$ do not intersect $N(A)$, 
consequently $y_p\in F_i'\cap F_j'$ and 
$y_p\notin F'$ because any facet of
$\Gamma'$ contains exactly one vertex of the edge $\{x_p,y_p\}$. Thus
$y_p=v$  and $v\in F_i'$, a contradiction. This proves $F'\cap
N(A)=\emptyset$, as required. Thus by reordering, we have that
$F_1',\ldots,F_r'$ is a 
shelling for the simplicial complex they generate. 
It is rapidly seen that $F_1,\ldots,F_r$ is also a shelling for the
simplicial complex they generate.

Next we consider the graph $G''=G\setminus\{x_g,y_g\}$ and the 
complex $\Gamma''=\Delta_{G''}^{[g-1]}$. Let
$F_1'',\ldots,F_m''$ be the facets of 
$\Gamma''$. By the induction hypothesis $\Gamma''$ is shellable. Thus we
may assume that $F_1'',\ldots,F_m''$ is a shelling of $\Gamma''$. It
is not hard to see that  
$$
H_1=F_1''\cup\{x_g\},\ldots,H_m=F_m''\cup\{x_g\}
$$
is the set of facets of $\Gamma$ containing $x_g$, and that
$H_1,\ldots,H_m$ is a shelling of the simplicial complex generated by
them. To finish the proof notice that 
$$
H_1,H_2,\ldots,H_m,F_1,F_2,\ldots,F_r,
$$
is clearly the complete list of facets of $\Gamma$ and they form a
shelling of $\Gamma$. Indeed for any $F_j$ one has that
$H_k=(F_j\setminus\{y_g\})\cup\{x_g\}$ is a facet of $\Gamma$ 
with  $F_i\setminus H_k=\{y_g\}$ and $H_k<F_j$.

$\Rightarrow$) The proof is by induction on $g$. We claim that $G$ has a
vertex of degree $1$. Let $F_1,\ldots,F_s$ be a shelling of $\Gamma$.
As $\{y_1,\ldots,y_g\}$ and $\{x_1,\ldots,x_g\}$ are facets of
$\Gamma$, we may
assume that $F_i=\{y_1,\ldots,y_g\}$, $F_j=\{x_1,\ldots,x_g\}$ and
$i<j$. Then there is $x_k\in F_j\setminus F_i$ and $F_\ell$ with
$\ell\leq j-1$ such that $F_j\setminus F_\ell=\{x_k\}$. 
Then 
$$\{x_1,\ldots,x_{k-1},x_{k+1}\ldots,x_g\}\subset F_\ell$$
and there is $y_t$ in $F_\ell$ for some $1\leq t\leq g$. Since
$$F_\ell=\{x_1,\ldots,x_{k-1},y_t,x_{k+1},\ldots,x_g\}$$ 
is an independent set of $G$, we get that
$y_t$ can only be adjacent to $x_t$. Thus $\deg(y_t)=1$ 
because $G$ has no isolated vertices. Thus we may order
$e_1,\ldots,e_g$ so that $\deg(x_g)=1$. Consider the graph 
$G'=G\setminus\{x_g,y_g\}$. Using \cite[Theorem~2.9]{bipartite-scm} 
we obtain that $\Delta_{G'}^{[g-1]}$ is a shellable complex. Hence by
induction hypothesis we can order $e_1,\ldots,e_{g-1}$ so that if
$\{x_i,y_j\}\in E(G')$, then $1\leq i\leq j\leq g-1$. To finish the 
proof note that any edge of $G$ is either an edge of $G'$ or 
an edge of $G$ containing $y_g$. \end{proof}

Some characterizations of condition 
{\rm ($\mathrm{h}_1$)} have been shown by Yassemi (personal
communication), and by  
Carr\`a Ferro and Ferrarello \cite{carra-ferrarello}. In
\cite{bipartite-scm} it is 
shown that if $G$ has a perfect matching and $R/I(G)$ is sequentially
Cohen-Macaulay, then condition ($\mathrm{h}_1$) holds. 

\begin{example}[\cite{bipartite-scm}]\label{bipscm-example} Let $G$
be the following bipartite graph. The ring 
$R/I(G)$ is not sequentially Cohen-Macaulay \cite{bipartite-scm} but 
the complex $\Delta_G^{[5]}$ is shellable.  

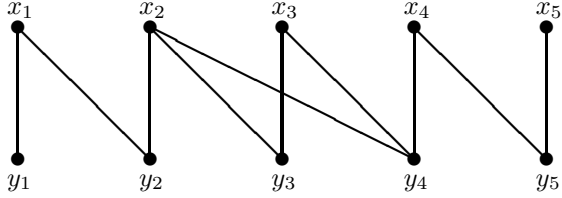
\begin{small}
\begin{picture}(100,60)(-85,20)
\thicklines
\put(0,0){\line(0,1){50}}
\put(0,0){\circle*{5}}
\put(0,50){\circle*{5}}
\put(0,50){\line(1,-1){50}}
\put(50,50){\line(1,-1){50}}
\put(50,50){\line(2,-1){100}}
\put(100,50){\line(1,-1){50}}
\put(50,0){\line(0,1){50}}
\put(100,0){\line(0,1){50}}
\put(150,0){\line(0,1){50}}
\put(200,0){\line(0,1){50}}
\put(150,50){\line(1,-1){50}}
 
\put(50,0){\circle*{5}}
\put(100,0){\circle*{5}}
\put(150,0){\circle*{5}}
\put(200,0){\circle*{5}}
\put(50,50){\circle*{5}}
\put(100,50){\circle*{5}}
\put(150,50){\circle*{5}}
\put(200,50){\circle*{5}}

\put(-4,55){$x_1$}
\put(46,55){$x_2$}
\put(96,55){$x_3$}
\put(146,55){$x_4$}
\put(196,55){$x_5$}

\put(-4,-10){$y_1$}
\put(46,-10){$y_2$}
\put(96,-10){$y_3$}
\put(146,-10){$y_4$}
\put(196,-10){$y_5$}
\end{picture}
\end{small}
\vspace{1.7cm}

\noindent A shelling of the facets of $\Delta_G^{[5]}$ is: 
$$
\begin{array}{ccc}
\{x_1, x_2, x_3, x_4, x_5\}<&\{x_2,x_3,x_4,x_5,y_1\}< &
\ \ \{x_3,x_4,x_5,y_1,y_2\}<\\ 
\{x_4,x_5,y_1,y_2,y_3\}<
&\{x_5,y_1,y_2,y_3,y_4\}<&\{y_1,y_2,y_3,y_4,y_5\} . 
\end{array}
$$
\end{example}

\begin{corollary}\label{herzog-hibi-ref} $G$ is a 
Cohen-Macaulay bipartite graph if and only if: 
{\rm ($\mathrm{h}_1'$)} $\Delta_G^{[g]}$ is shellable, 
$g={\rm ht}\, I(G)$, 
and {\rm ($\mathrm{h}_2'$)} $G$ is unmixed.
\end{corollary} 

\begin{proof}
It follows using Corollary~\ref{unmixed-bip} together with
Theorems~\ref{herzog-hibi} and \ref{ordering-char}.  
\end{proof}

This corollary shows that $G$ is Cohen-Macaulay if and only if
$\Delta_G$ is pure shellable \cite{EV,bipartite-scm}.

The natural generalization of a bipartite graph is a balanced 
clutter. The next example shows that Theorems~\ref{ordering-char} and 
\ref{herzog-hibi} do not extend to balanced
clutters.

\begin{example}\label{counterexample} 
Consider the clutter $\mathcal{C}$ whose edge ideal is generated by:
$$
\begin{array}{ccc}
a_1b_1c_1d_1g_1h_1k_1,&a_2b_2c_2d_2g_2h_2k_2,&a_3b_3c_3d_3g_3h_3k_3,\\
a_4b_4c_4d_4g_4h_4k_4,&a_1b_1c_1d_1g_2h_3k_4,&a_1b_2c_3d_4g_2h_3k_4.
\end{array}
$$
This clutter is balanced. Indeed its incidence matrix $A$ is totally
unimodular, i.e., each $i\times i$ minor of $A$ is $0$ or $\pm 1$ for
all $i\geq 1$. Furthermore $\mathcal{C}$ satisfies condition 
(b) of Corollary~\ref{char-tbc}.
Hence $I(\mathcal{C})$ is Cohen-Macaulay. However we cannot order 
its vertices so that it becomes 
an admissible uniform clutter.
\end{example}
\noindent
{\bf Acknowledgments.} We gratefully acknowledge the computer algebra 
system {\tt Macaulay\/}$2$ \cite{mac2} which was invaluable in our
work on this paper. 
The second author also acknowledges the financial support of
COFAA-IPN and SNI; the third author acknowledges the financial
support of  
CONACyT grant 49251-F and SNI. We also thank the referee for a
careful reading of the paper and for the  improvements suggested.

\bibliographystyle{plain}

\end{document}